\documentclass[twoside, 10pt]{article}

\usepackage{mathrsfs,amsfonts,amsmath}
\RequirePackage{ifthen,calc}
\usepackage{theorem}
\usepackage[dvips]{color}
 \catcode`@=11
 \def\@evenhead{\hbox to\textwidth{\footnotesize\rm\thepage \hfill
  {\it Li  and  Xiao}}} 

 \def\@oddhead{\hbox to \textwidth{\footnotesize{\it
   Fluctuation Limits of Weakly Degenerate Branching Systems } \hfill\thepage}}

 \renewcommand{\section}{\makeatletter
 \renewcommand{\@seccntformat}[1]{{\csname the##1\endcsname.}\hspace{0.45em}}
 \makeatother \@startsection
{section}
{1}
{0pt}
{\baselineskip}
{0.5\baselineskip}
{\normalsize\bfseries\mathversion{bold}}}
\catcode`@=12
\newcommand\ack{\section*{Acknowledgement}}

\newtheorem{thm}{\noindent Theorem}[section]
\newtheorem{lem}{\noindent Lemma}[section]

\newtheorem{prop}{\noindent Proposition}[section]

\theoremheaderfont{\normalfont\bfseries} \theorembodyfont{\slshape}
{\theorembodyfont{\rmfamily}} {\theorembodyfont{\rmfamily}
}
{\theorembodyfont{\rmfamily} } {\theorembodyfont{\rmfamily}
\newtheorem{rem}{\noindent Remark}[section]}
{\theorembodyfont{\rmfamily} } {\theorembodyfont{\rmfamily} } {\theorembodyfont{\rmfamily}}
 \def\beqlb{\begin{eqnarray}}\def\eeqlb{\end{eqnarray}}
 \def\beqnn{\begin{eqnarray*}}\def\eeqnn{\end{eqnarray*}}
 \numberwithin{equation}{section}
 \def\qed{\hfill$\square$\smallskip}
 \def\R{{\mathbb R}}
 \def\bfE{{\mathbb{E}}}
 \def\e{\mathrm{e}}
 
 \def\d{\mathrm{d}}

\begin{document}

\renewcommand{\baselinestretch}{1.3}
\normalsize

\title{\LARGE\bf Occupation Time Fluctuations of
Weakly\\ Degenerate Branching Systems }
\author{Yuqiang Li\thanks{Research partially
        supported by a grant from NSF of China 10901054.}
\\ \small School of Finance and Statistics, East China Normal University
\\ \small Shanghai 200241, P. R. China.\\
\\  Yimin Xiao\thanks{Research partially
        supported by  NSF grant DMS-1006903.}
\\ \small Department of Statistics and Probability, Michigan State
University
\\ \small East Lansing, MI 48824, USA}

\date{}

\maketitle

\centerline{\textbf{Abstract}}

\noindent{We establish limit theorems for re-scaled occupation time fluctuations
of a sequence of branching particle systems in $\R^d$ with anisotropic
space motion and weakly degenerate splitting ability. In the case of
large dimensions, our limit processes lead
to a new class of operator-scaling Gaussian random fields with
non-stationary increments. In the intermediate and critical dimensions,
the limit processes have spatial structures analogous to (but
more complicated than) those arising from the
critical branching particle system without degeneration considered
by Bojdecki et al. \cite{BGT061,BGT062}. Due to the weakly degenerate
branching ability, temporal structures of the limit processes in all three
cases are different from those obtained by Bojdecki et al.
\cite{BGT061,BGT062}.}

\bigskip


\bigskip


\medskip

\newpage

\renewcommand{\baselinestretch}{1.0}
\normalsize

\section{Introduction}
Consider a kind of branching particle systems in $\R^d$ as follows. Particles
start off at time $t=0$ from a Poisson random field with Lebesgue
intensity measure $\lambda$ and evolve independently. They move in $\R^d$
according to a L\'evy process
$$
\vec{\xi}=\{\vec{\xi}(t),
t\geq0\}=\{(\xi_1(t),\xi_2(t),\cdots,\xi_d(t)), t\geq 0\}
$$
with independent stable components as in \cite{PT69}, i.e. for
every $1 \le k\leq d$, $\xi_k=\{\xi_k(t), t\geq0\}$ is a real-valued symmetric
$\alpha_k$-stable L\'{e}vy process, and $\xi_1,\cdots,\xi_d$ are
independent of each other. In addition, the particles split at a rate
$\gamma$ and the branching law at age $t$ has the generating function
$$
g(s,t)=\Big(1-\frac{\e^{-\delta t}}{2}\Big)
 +\e^{-\delta t}\frac{s^2}{2}, \qquad 0\leq s\leq1,\  t\geq 0.
$$
Intuitively, in this model, the particles' motion in $\R^d$ is anisotropic (i.e.,
the motion in different direction is controlled by different mechanism) and their
ability of splitting new particles declines as their ages increase (at time $t$
each particle produces 2 particles with probability $\e^{-\delta t}/2$ and no particles
with probability $1 -\e^{-\delta t}/2$). For
simplicity of notation, we denote the vector $(\alpha_1,\cdots,\alpha_d)$ by
$\vec{\alpha}$, and call this model a $(d, \vec{\alpha},
\delta,\gamma)$-degenerate branching particle system.

Let $N(s)$ denote the empirical measure of the particle system at
time $s$, i.e. $N(s)(A)$ is the number of particles in the set
$A\subset\R^d$ at time s.  We call the measure-valued process
$$
L(t)=\int_0^{t}N(s)\d s,\quad t\geq0,
$$
the occupation time of the system, and call the process
$$
X(t)=\int_0^{t}(N(s)-\bfE(N(s)))\d s,\quad t\geq0,
$$
the occupation time fluctuation, where $\bfE(N(s))$ is the
expectation functional understood as $\langle
\bfE(N(s)),\phi\rangle=\bfE(\langle N(s), \phi\rangle)$ for all
$\phi\in\mathcal{S}(\R^d)$, the space of smooth rapidly decreasing
functions. Here and sometimes in the sequel, we write $\langle\mu,
f\rangle=\int f\d\mu$, where $\mu$ is a measure and $f$ a
measurable function.

The main purpose of this paper is to study the scaling limit of
occupation time fluctuations of a sequence of $(d, \vec{\alpha},
\delta,\gamma)$-degenerate branching particle systems (the specific
assumptions will be given below). This is motivated by the limit theorems
for occupation time fluctuations of branching particle
systems established recently by Bojdecki et al. \cite{BGT061,BGT062,BGT071,BGT072}
and by the research on anisotropic Gaussian random fields (cf. \cite{BMS07,BE03,X09}).

Branching particle systems and their associated
superprocess have been widely studied; see for example, \cite{DP99,
GW91,I86,Le99,L98,M06}. 
Recently Bojdecki et al. \cite{BGT04,BGT061,BGT062,
BGT071,BGT072,BGT081, BGT082} established very interesting limit theorems for
occupation time fluctuations of the $(d, \alpha,\beta)$-branching systems,
and showed that the limit processes depend
on the relations among the parameters $d, \alpha$ and $\beta$ as well as
the initial distribution of the particles. Under the assumptions that the
initial state $N(0)$ is a Poisson random measure with intensity
the Lebesgue measure $\lambda$, particles move independently following
isotropic $\alpha$-stable L\'evy processes in $\R^d$ and undergo critical
binary branching (that is, $\beta = 1$), Bojdecki et al. \cite{BGT061, BGT062}
proved that, if the dimension is intermediate (i.e., $\alpha<d<2\alpha$),
then the limit process has the form
$C\lambda \xi$, where $C$ is a constant and $\xi=(\xi_t,t\geq0)$ is
a sub-fractional Brownian motion (sub-fBm) defined in \cite{BGT04};
if the dimension is critical (i.e., $d = 2\alpha$), then the limit process
has the form $C\lambda \zeta$, where $C$ is a constant and $\zeta=\{\zeta_t,
t\geq0\}$ is a standard Brownian motion; and if the dimension is large
(i.e., $d > 2 \alpha$), then the limit is a generalized Gaussian process
which is essentially a generalized Wiener process whose
spatial structure is determined solely by $\alpha$.

We observe that the limit process, say $X$, of occupation time
fluctuation is a process valued in $\mathcal{S}'(\R^d)$, the dual
space of $\mathcal{S}(\R^d)$. Intuitively, when we consider a multi-parameter
process $Y=\{Y(z), z\in\R_+^d\}$ defined by $Y_z=\langle X(1), {\bf
1}_{D(z)}\rangle$, where $D(z)=\{(x_1,\cdots,x_d), 0\leq x_k\leq
z_k, 1\leq k\leq d\}$, we obtain a real-valued random field.
For random fields, ``anisotropy" is a distinct property from those of one-parameter
processes and is important for many applied areas such as geophysical, economic and
ecological sciences. Several classes of anisotropic Gaussian random
fields have been explicitly constructed by using stochastic integrals and their
properties have been studied. We refer to \cite{BMB06,BMS07,BE03,LX09,X09,XZ04}
and the references therein for further information. Many of the anisotropic
random fields $Y = \{Y(t),\, t \in \R^d\}$ in the literature have the following
scaling property: There exists a linear operator $E$ (which may not be unique)
on $\R^d$ with positive real parts of the eigenvalues such that for all constants $c > 0$,
\begin{equation}\label{def:OSS}
\big\{ Y(c^E\,t),\, t \in \R^d\big\} \stackrel{f.d}{=}
\big\{c\,Y(t),\, t \in \R^d\big\}.
\end{equation}
Here and in the sequel, ``$\stackrel{f.d}{=}$" means equality in all
finite-dimensional distributions and, for $c > 0$, $c^E$ is the linear
operator on $\R^d$ defined by $c^E = \sum_{n=0}^\infty \frac{(\ln c)^n E^n} {n!}$.
If $Y = \{Y(t),\, t \in \R^d\}$ satisfies (\ref{def:OSS}), we call $Y$ an operator-scaling
random field with (time-scaling) exponent $E$. More general forms of scaling properties
are also possible for multivariate random fields (cf. \cite{LX09}).

It is known that some anisotropic Gaussian random fields such as the fractional Brownian
sheets can be obtained as scaling limit of partial sums of discrete-time random fields
(\cite{La07}). Such results provide physical interpretation for anisotropic random fields.
Given the significance of branching particle systems, it is interesting to
investigate whether operator-scaling random fields arise naturally in such models.
Considering branching particle systems with anisotropic particle motions is a natural
step in this direction.

Now we specify our setting. Because of the sub-critical branching laws at positive ages,
a fixed $(d, \vec{\alpha}, \delta,\gamma)$-branching particle system with $\delta>0$
will go to local extinction as time elapses. To overcome this difficulty, we borrow
the idea of nearly critical branching processes (see \cite{IPZ05,Li09,ST94}) and consider
a sequence of $(d,\vec{\alpha}, \delta_n,\gamma)$-branching particle systems with
$\delta_n\to0$ as $n\to\infty$. We study the limit process of the re-scaled occupation
time fluctuations
$$
X_n(t)=\frac{1}{F_n}\int_0^{nt}(N_n(s)-\bfE(N_n(s)))\d s,
$$
where $\{F_n, n \ge 1\}$ is a sequence of norming constants to be chosen appropriately.
This setting is different from those in Bojdecki et al. \cite{BGT061, BGT062} and
other aforementioned references, where a fixed system is studied.
We further assume that $n\delta_n\to\theta$ for some constant $\theta\geq0$ as $n\to\infty$,
which is referred to as weak degeneration. Our focus is on how the anisotropy of the space
motion $\vec{\xi}$ and the degeneration of splitting ability affect the limit
process of $X_n$. Let
$$\bar{\alpha}=\sum_{k=1}^d\frac{1}{\alpha_k}.$$
Our results will show that the scaling limit of $X_n(t)$ depends crucially on
$\bar{\alpha}$. We assume throughout the paper that $\bar{\alpha}>1$. In analogy
to the terminology of Bojdecki et al. \cite{BGT061, BGT062}, we refer to the three cases
$\bar{\alpha} >2$, $\bar{\alpha} =2$ and $1 <\bar{\alpha} <2$ as the large dimension,
critical dimension and intermediate dimension, respectively.

The remainder of this paper is organized as follows. In Section 2,
we derive an explicit expression for the re-scaled occupation time fluctuation
$X_n$ and state our main theorems. These results show that operator-scaling
Gaussian random fields can arise only in the large dimension case and they
generally have non-stationary increments. In the cases $\bar{\alpha}\in(1, 2)$
and $\bar{\alpha}=2$, the spatial structures of the limit processes are similar to (but more
complicated than) those of the $(d, \alpha, 1)$-branching systems in \cite{BGT061,BGT062}.
However, the temporal structures of limit processes are always different. See Remark
\ref{s2-rm-2} for more details.

The space-time method employed in this paper for proving our main results
are analogous to those developed by Bojdecki et al. \cite{BGT061, BGT062}, with some
nontrivial modifications to handle the new technical complexities caused
by the anisotropy of the space motion and the degenerate branching
ability. Section 3 contains discussions on the Laplace functionals of
the occupation time fluctuations and is devoted to the proofs of several
technical lemmas. Finally in Section 4 we provide proofs of the
main results stated in Section 2.

Unless stated otherwise, $K$ denotes an unspecified positive finite
constant which may not necessarily be the same in each occurrence.

\section{Main results}

Consider a sequence of $(d,
\vec{\alpha},\delta_n,\gamma)$-degenerate branching particle
systems. For every $n \ge 1$, let $N_n(t)$ be the corresponding
empirical measures and let $H$ be the $d\times d$ diagonal matrix
$(1/\alpha_k)_{1\leq k\leq d}$. The corresponding space motion is
denoted by $\vec{\xi}_n = \{\vec{\xi}_n(t), t \ge 0\}$. We assume
that $\{\vec{\xi}_n, n \ge 1\}$ is a sequence of identically
distributed $\R^d$-valued L\'evy processes with $\alpha_k$-stable
components ($1\le k \le d$). The distribution of $\vec{\xi}_n$ is
completely determined by its characteristic function
\begin{equation}\label{Eq:ChF}
{\mathbb E}\Big(\e^{i \langle z, \vec{\xi}_n(t)\rangle}\Big) =
\e^{- t \sum_{k=1}^d |z_k|^{\alpha_k}}, \quad z \in \R^d.
\end{equation}
It follows that $\vec{\xi}_n$ has the following operator-self-similarity:
For all constants $c > 0$,
\[
\{\vec{\xi}_n(ct), t \ge 0\} \stackrel{f.d}{=} \{c^H\vec{\xi}_n(t), t \ge 0\}.
\]
Hence, $\vec{\xi}_n$ is an operator-stable L\'evy process on $\R^d$ (cf.
\cite[Theorem 7]{HM82} or \cite{Sato}). We refer to \cite{Sato} for a
systematic account on L\'evy processes. In the following, we denote the
semigroup of $\vec{\xi}_n$ by $\{T_t\}_{t\geq0}$, i.e.,
$$ T_sf(x):=\bfE(f(\vec{\xi}_n(t+s))|\vec{\xi}_n(t)=x)$$
for all $s, t\geq0$, $x\in\R^d$ and bounded measurable functions
$f$ on $\R^d$. To avoid misunderstanding, we sometimes write $T_sf(x)$ as
$T_s(f(\cdot))(x)$.

For convenience, we use the notation
$$\bfE_x(f(N_n))=\bfE(f(N_n)|N_n(0)=\epsilon_x)$$
for any measurable function $f$ on $D([0,\infty),
\mathcal{S}'(\R^d))$, where $\epsilon_x$ denotes the unit measure
concentrated at $x\in\R^d$. For all $\phi \in \mathcal{S}(\R^d)$, define
\beqlb\label{s2-1}
 F_{n,\phi}(x,s)=\bfE_x(\langle N_n(s), \phi\rangle).
 \eeqlb
Note that all particles split at the rate $\gamma$ and evolve
independently. By the renewal method, it is not hard to verify that
for all $\phi \in \mathcal{S}(\R^d)$,
 \beqnn
  F_{n,\phi}(x, s)&=&\e^{-\gamma s}T_s\phi(x)+\int_0^s\gamma\e^{-\gamma
 u}\e^{-\delta_n u}T_uF_{n,\phi}(\cdot, s-u)(x)\d u\nonumber\\
 &=&\e^{-\gamma s}T_s\phi(x)+\int_0^s\gamma\e^{-(\gamma+ \delta_n)(s-u)}
 \,T_{s-u}F_{n,\phi}(\cdot, u)(x)\d  u.
 \eeqnn
Therefore,
 \beqlb\label{s2-2}
 \frac{\partial F_{n,\phi}}{\partial
 s}&=&(-\gamma+A)F_{n,\phi}+\gamma
 F_{n,\phi}-\delta_n(F_{n,\phi}-\e^{-\gamma s}T_s\phi(x))\nonumber
 \\&=&(A-\delta_n)F_{n,\phi}+\delta_n\e^{-\gamma s}T_s\phi(x),
 \eeqlb
where $A$ denotes the infinitesimal generator of $T$. Note that
$F_{n,\phi}(x, 0)=\phi(x)$. From (\ref{s2-2}), it follows that
 \beqlb \label{s2-3}
 F_{n, \phi}(x, s)&=&\e^{-\delta_n
 s}T_s\phi(x)+\delta_n\int_0^s\e^{-\delta_n(s-u)}T_{s-u}(\e^{-\gamma
 u}T_u\phi)(x)\d u\nonumber
 \\&=&\e^{-\delta_n
 s}T_s\phi(x)\Big(1+\frac{\delta_n}{\gamma-\delta_n}(1-\e^{-(\gamma-\delta_n)s})\Big).
 \eeqlb
Let
 \beqlb\label{s2-4}
 \bar{f}_n(s):=1+\frac{\delta_n}{\gamma-\delta_n}(1-\e^{-(\gamma-\delta_n)s})\ \
 \text{ and }\;\ \  f_n(s)=\bar{f}_n(s)\e^{-\delta_n s}.
 \eeqlb
Then, for all $n\geq1, x\in\R^d$ and $\phi\in\mathcal{S}(\R^d)$,
(\ref{s2-3}) can be written as \beqlb\label{s2-5}
 F_{n,\phi}(x, s)=f_n(s)T_s\phi(x).
\eeqlb
Because the Lebesgue measure is an invariant measure for symmetric
stable L\'{e}vy processes and all components of $\vec{\xi}_n$ are
symmetric stable L\'{e}vy processes and independent of each other,
we use Fubini's theorem to derive that
 \beqlb\label{s2-6}
 \int_{\R^d} T_s\phi(x)\d x=\int_{\R^d} \phi(x)\d x
 \eeqlb
for all $s>0$ and $\phi\in\mathcal{S}(\R^d)$. Therefore, from the
fact that $N_n(0)$ has a Poisson distribution with Lebesgue
intensity measure $\lambda$, we find that
 $$
 \bfE(\langle N_n(s), \phi\rangle)=\int_{\R^d} F_{n, \phi}(s, x)\d
 x=f_n(s)\int_{\R^d} \phi(x)\d x=f_n(s)\langle \lambda, \phi\rangle.
 $$
Now we define the occupation time fluctuation process $X_n=\{X_n(t),
t\geq0\}$ as follows
 \beqlb\label{s2-7}
 \langle X_n(t),\phi\rangle=\frac{1}{F_n}\int_0^{nt}\langle N_n(s)-f_n(s)\lambda,
 \phi\rangle\d s
 \eeqlb
for every $\phi\in\mathcal{S}(\R^d)$, where $F_n$ is a suitable
scaling constant.

Below, we always assume that $N_n(s)$ is the empirical measure of a $(d, \vec{\alpha},
\delta_n,\gamma)$-degenerate branching particle system with $\bar{\alpha}=
\sum_{i=1}^d\alpha_k^{-1}>1$ and there is a constant $\theta\in[0,\infty)$ such that
$n\delta_n\to\theta$ as $n\to\infty$.
Let $\widehat{\phi}(z)$ $(z\in\R^d)$ denote the Fourier
transform of function $\phi\in L^1(\R^d)$, i.e.,
 $\widehat{\phi}(z)=\int_{\R^d}\e^{i\langle x, z\rangle}\phi(x)\d x$.

We distinguish three cases: (i) $\bar{\alpha}>2$ (the large dimension case);
(ii)  $\bar{\alpha}=2$ (the critical dimension case) and (iii)  $1 <\bar{\alpha}<2$
(the intermediate dimension case).
The main results of this paper are stated as follows.

\begin{thm}\label{s2-thm-1}
When $\bar{\alpha}>2$, let  $F_n^2=n$. Then for every $t>0$ and
$\psi\in\mathcal{S}(\R^{d+1})$,
$\int_0^t\langle X_n, \psi(\cdot, s)\rangle\d s$ converges in distribution
to $\int_0^t\langle X, \psi(\cdot, s)\rangle\d s$,
where $X$ is a centered Gaussian process with covariance function
\beqlb\label{s2-thm-1-1}
&&{\rm Cov} (\langle X(r),\phi_1\rangle,\langle X(t), \phi_2\rangle)\nonumber \\
&&\quad=\frac{1}{(2\pi)^{d}}\int_{\R^d}\bigg(\frac{2}{\sum_{k=1}^d|z_k|^{\alpha_k}}
 +\frac{\gamma}{(\sum_{k=1}^d|z_k|^{\alpha_k})^2}\bigg)\widehat{\phi}_1(z)
 \overline{\widehat{\phi}_2(z)}\d z\int_0^{r\wedge t}\e^{-\theta s}\d s.
\eeqlb
\end{thm}

Theorem \ref{s2-thm-1} shows that the limit process $X$ can be decomposed as
the sum of two independent $\mathcal{S}'(\R^d)$-valued Gaussian
processes, say $\bar{X}_1$ and  $\bar{X}_2$, where
\beqnn
 {\rm Cov} (\langle \bar{X}_1(r),\phi_1\rangle,\langle \bar{X}_1(t), \phi_2\rangle)=\frac{1}{(2\pi)^{d}}\int_{\R^d}\frac{2}{\sum_{k=1}^d|z_k|^{\alpha_k}}
 \widehat{\phi}_1(z)\overline{\widehat{\phi}_2(z)}\d
 z\int_0^{r\wedge t}\e^{-\theta s}\d s
\eeqnn
 and
\beqnn
  {\rm Cov} (\langle \bar{X}_2(r),\phi_1\rangle,\langle \bar{X}_2(t), \phi_2\rangle)=\frac{1}{(2\pi)^{d}}\int_{\R^d}\frac{\gamma}{(\sum_{k=1}^d|z_k|^{\alpha_k})^2}
 \widehat{\phi}_1(z)\overline{\widehat{\phi}_2(z)}\d
 z\int_0^{r\wedge t}\e^{-\theta s}\d s.
 \eeqnn
Due to the nuclear property of $\mathcal{S}(\R^d)$, we can define $Y_1(u)=\langle \bar{X}_1(1), {\bf
1}_{D(u)}\rangle$ and $Y_2(u)=\langle \bar{X}_2(1), {\bf 1}_{D(u)}\rangle$
for all $u=(u_k)_{1\leq k\leq d}\in[0,\infty)^d$, where
$D(u)=\{(x_1,\cdots,x_d), 0\leq x_k\leq u_k, 1\leq k\leq d\}$. The following proposition characterize
the operator-scaling properties of the Gaussian random fields $Y_1$ and $Y_2$.

\begin{prop}\label{s2-prop} Suppose $\bar{\alpha}>2$.
\noindent(1) When $d=1$, denote $\vec{\alpha}$ by $\alpha$. Then
$Y_1=\{Y_1(u), u\geq0\}$ and $Y_2=\{Y_2(u), u\geq0\}$, up to
some multiplicative constants, are the fractional Brownian motion
with Hurst indices $H_1=\frac{1+\alpha}{2}$ and $H_2=\frac{1+2\alpha}{2}$, respectively.

\noindent(2) When $d\geq2$, $Y_1=\{Y_1(u), u\in[0,\infty)^d\}$ and $Y_2=\{Y_2(u),
u\in[0,\infty)^d\}$ are operator-scaling Gaussian random fields with exponent
$H$ (which is the $d\times d$ diagonal matrix $(1/\alpha_k)_{1\leq
k\leq d}$) and have non-stationary increments.
\end{prop}

From the proof of Proposition \ref{s2-prop} in Section 4, one can see that,
if $\alpha_1, \ldots, \alpha_d$ are not the same, then the Gaussian random fields
$Y_1$ and $Y_2$ are anisotropic. However, they are different from the fractional
Brownian sheets or the operator-scaling Gaussian fields constructed by Bierm\'e et al.
\cite{BMS07}, and can be further studied by applying the general methods in Xiao
\cite{X09}.

\begin{thm}\label{s2-thm-2}
When $\bar{\alpha}=2$, let  $F_n^2=n\ln n$. Then for any $t>0$ and $\psi\in
\mathcal{S}(\R^{d+1})$,
$\int_0^t\langle X_n, \psi(\cdot, s)\rangle\d s$ converges in
distribution to $\int_0^t\langle X, \psi(\cdot, s)\rangle\d s$,
where
$$ X=\sqrt{\frac{\gamma}{(2\pi)^{d}}}\ \lambda\xi_\theta.
$$
In the above
$\xi_{\theta}=\{\xi_\theta(t), t\geq 0\}$ is a centered Gaussian
process with covariance function
\beqlb\label{s2-thm-2-1}
C_\theta(r, t)=\bfE(\xi_\theta(r)\xi_\theta(t))=\int_0^{r\wedge t} \e^{-\theta s}\d
s\int_{\R^d}\varphi(\theta, r-s, t-s, \sum_{k=1}^d|y_k|^{\alpha_k})\d y,
\eeqlb
where for any $x, u, v\geq 0$ and $z>0$,
\beqlb\label{def1}
\varphi(x, u, v, z)&=&\frac{uz\e^{-u(x+z)}(1-\e^{-v(x+z)})}
{(x+z)^2}+\frac{vz\e^{-v(x+z)} (1-\e^{-u(x+z)})}
 {(x+z)^2}\nonumber
 \\&&+\frac{2x(1-\e^{-u(x+z)})
 (1-\e^{-v(x+z)})}{(x+z)^3}.
 \eeqlb
\end{thm}

\begin{thm}\label{s2-thm-3}
When $1<\bar{\alpha}<2$, let $F_n=n^{(3-\bar{\alpha})/2}$. Then
$X_n\Rightarrow X$ in $C([0, 1], \mathcal{S}'(\R^d))$ as $n\to\infty$,
where $X=\{X(t), t\in[0,
1]\}$ is a centered Gaussian process with covariance function
$$
{\rm Cov}(\langle X(r),\phi_1\rangle,\langle X(t), \phi_2\rangle)=
 \frac{\gamma}{\pi^{d}}\prod_{k=1}^d\frac{\Gamma(1/\alpha_k)}{\alpha_k}
 C(r, t)\langle\lambda, \phi_1\rangle\langle\lambda,\phi_2\rangle
$$
with
\beqlb\label{s2-thm-3-2}
C(r, t)=\int_0^r\e^{-\theta u}\d u\int_0^t\e^{-\theta v}\d
v\int_0^{u\wedge v} \frac{\e^{\theta s}}{(u+v-2s)^{\bar{\alpha}}}\d
s.
\eeqlb
\end{thm}

\begin{rem}
According to Definition 2.1 in Bojdecki et al. \cite{BGT072}, the convergence in Theorem
\ref{s2-thm-1} and Theorem \ref{s2-thm-2} is called the convergence in the integral
sense. We point out that it is possible to prove the tightness of the processes
$X_n$ in $C([0, 1], \mathcal{S}'(\R^d))$ by using the method in \cite{BGT062}
and hence to strengthen the conclusions of Theorem \ref{s2-thm-1} and  \ref{s2-thm-2}
to weak functional convergence in $C([0, 1],\mathcal{S}'(\R^d))$. Since this
approach is lengthy and tedious, in order to save the space of this paper,
we do not prove this stronger sense of convergence and focus on identifying
the limit processes in Theorems \ref{s2-thm-1} and \ref{s2-thm-2}.
\end{rem}

Comparing our results with those of Bojdecki et al. \cite{BGT061,BGT062}
on the critical binary branching particle systems in $\R^d$ with symmetric
stable L\'evy motion, we have the following comments.

\begin{rem}\label{s2-rm-2}
(1) Let $\bar{X}$ denote a $\mathcal{S}'(\R^d)$-valued Gaussian process with
covariance function
\beqnn
&&{\rm Cov}(\langle \bar{X}(r),\phi_1\rangle,\langle \bar{X}(t), \phi_2\rangle)
=(r\wedge t)\int_{\R^d}\bigg(\frac{2}{\sum_{k=1}^d|z_k|^{\alpha_k}}
+\frac{\gamma}{(\sum_{k=1}^d|z_k|^{\alpha_k})^2}\bigg)\widehat{\phi}_1(z)
\overline{\widehat{\phi}_2(z)}\d z.
\eeqnn
Then $\bar{X}$ is an analogue of the limit process in \cite[Theorem 2.2 (a)]{BGT062}.
The limit process in Theorem \ref{s2-thm-1} can be written as
$$X(t)=\int_0^t\e^{-\frac{\theta}{2}s}\d  \bar{X}(s),$$
namely, for any $\phi\in \mathcal{S}(\R^d)$, $\langle X(t),
\phi\rangle=\int_0^t\e^{-\frac{\theta}{2}s}\d
\langle \bar{X}(s), \phi\rangle.$

(2) For $x, u, v\geq 0$ and $z>0$, let
$$\varphi_1(x, u, z)=\frac{uz\e^{-u(x+z)}}{(x+z)^2}\ \;\text{ and}\;\;
\varphi_2(x,u,v,z)=\frac{2x(1-\e^{-u(x+z)})
 (1-\e^{-v(x+z)})}{(x+z)^3}.$$
Then
\beqlb\label{rm-2-1}
\varphi(x, u,v,z)=\varphi_1(x, u,z)+\varphi_1(x,v,z)-\varphi_1(x,u+v,z)+\varphi_2(x,u,v,z).
\eeqlb
For $x\geq0$, define
$$\Phi(x):=\int_{\R^d}\varphi_1\Big(x, 1,\sum_{k=1}^d|y_k|^{\alpha_k}\Big)\d y.$$
Then
$$
\Phi(x) \leq\int_{\R^d}\frac{\e^{-\sum_{k=1}^d|y_k|^{\alpha_k}}}
{\sum_{k=1}^d|y_k|^{\alpha_k}}\d y<\infty,
$$
where the finiteness follows from Remark \ref{s2-rm-1} below. Then
for $x\geq 0$ and $s<r$, using the substitutions $y'=(r-s)^{H}y$ and
the condition $\bar{\alpha}=2$, we obtain that
 \beqlb\label{rm-2-2}
 \int_{\R^d}\varphi_1\Big(x, r-s,\sum_{k=1}^d|y_k|^{\alpha_k}\Big)\d y
 =\int_{\R^d}\varphi_1\Big((r-s)x, 1,\sum_{k=1}^d|y'_k|^{\alpha_k}\Big)\d y'=\Phi((r-s)x).
 \eeqlb
Let
$$ \Phi_1(\theta, r,s, t)=\Phi((t-s)\theta)+\Phi((r-s)\theta)-\Phi((t+r-2s)\theta).$$
By (\ref{rm-2-1}) and (\ref{rm-2-2}) one can verify that for all $r, t, s\in[0,1]$
with $s<r\wedge t$,
\beqnn
\int_{\R^d}\varphi\Big(\theta, r-s, t-s, \sum_{k=1}^d|y_k|^{\alpha_k}\Big)\d
 y&=&\Phi_1(\theta, r, s, t)+\Phi_2(\theta, r, s, t),
\eeqnn
where
$$
\Phi_2(\theta, r, s, t)
=\int_{\R^d}\varphi_2\Big(\theta, r-s, t-s, \sum_{k=1}^d|y_k|^{\alpha_k}\Big)\d
 y<\infty.
 $$
So the Gaussian process $\xi_\theta$ in Theorem \ref{s2-thm-2} has
the covariance function
   $$C_\theta(r, t)=\int_0^{r\wedge t}\e^{-\theta s}\Phi_1(\theta, r, s, t)\d s
   +\int_0^{r\wedge t}\e^{-\theta s}\Phi_2(\theta, r, s,
  t)\d s.$$
 When $\theta=0$, $C_0(r, t)=\Phi(0)(r\wedge t)$. In this case, up to a
 multiplicative constant, the
 limit process $X$ is Brownian motion, which is the same as that in \cite[Theorem 2.2 (b)]{BGT062}.
 However, when $\theta\not=0$, the limit process has complicated covariance function
 and can be expressed as the sum of two independent $\mathcal{S}'(\R^d)$-valued Gaussian
 processes--the first one is an extension of that in \cite[Theorem 2.2 (b)]{BGT062}
 and the second is a new process.

(3) In Theorem \ref{s2-thm-3}, if $\theta=0$, then, up to a
multiplicative constant, the limit process is the same as that in
Theorem 2.2 of \cite{BGT061}, i.e. the limit process can be written
as $K\lambda\xi^h$, where $\xi^h$ is a sub-fractional Brownian motion with covariance
function
$$
s^{3-\bar{\alpha}}+t^{3-\bar{\alpha}}-\frac{1}{2}[(s+t)^{3-\bar{\alpha}}+|s-t|^{3-\bar{\alpha}}].
$$
This sub-fBm is self-similar with index $(3-\bar{\alpha})/2$. When
$\theta \ne 0$, the Gaussian process with covariance
(\ref{s2-thm-3-2}) in the limit process in Theorem \ref{s2-thm-3} is
not self-similar and is a new process.
\end{rem}

We end this section with some preliminary facts which will be useful for our proofs.
The proof of Lemma \ref{s2-lm-1} is elementary and thus is omitted.
\begin{lem}\label{s2-lm-1}
{\rm (1)} If $\bar{\alpha}>2$, then
$$\int_{[0, 1]^d}\frac{1}{\sum_{k=1}^d z_k^{\alpha_k}}\d z\leq
  \int_{[0, 1]^d}\frac{1}{\sum_{k=1}^d z_k^{2\alpha_k}}\d z<\infty.$$

{\rm (2)} If $\bar{\alpha}<2$, then
$$\int_{\R_+^d\setminus[0, 1]^d}\frac{1}{\sum_{k=1}^d z_k^{2\alpha_k}}\d z<\infty.$$
\end{lem}

\begin{rem}\label{s2-rm-1}
Lemma \ref{s2-lm-1} implies that for $r\in(0, \bar{\alpha})$,
 $ \int_{[-1, 1]^d}\frac{1}{\sum_{k=1}^d|z_k|^{r\alpha_k}}\d z<\infty,
 $
and for $r>\bar{\alpha}$,
 $ \int_{\R^d\setminus[-1,
1]^d}\frac{1}{\sum_{k=1}^d|z_k|^{r\alpha_k}}\d z<\infty.
 $
 Therefore, if $\tau(z)$
is bounded and $\int_{\R^d} \tau(z)\d z<\infty$, then
  $ \int_{\R^d}\frac{\tau(z)}{\sum_{k=1}^d|z_k|^{r\alpha_k}}\d
 z<\infty$
for all $r\in(0, \bar{\alpha})$.
\end{rem}

We will repeatedly use the following formulas
involving Fourier transforms. Let $\phi_1$, $\phi_2$ and $\phi_3$
be functions from $\R^d$ to $\R$, bounded and integrable. Then
\beqlb
 \int_{\R^d}\phi_1(x)\phi_2(x)\d x&=& \frac{1}{(2\pi)^{d}}\int_{\R^d}
 \widehat{\phi}_1(z)\overline{\widehat{\phi}_2(z)}\d
 z,\label{s2-12}
 \eeqlb
(the Plancherel formula). If $\widehat{\phi_1}$ and $\widehat{\phi_2}$
are integrable, then
\beqlb
\int_{\R^d}\phi_1(x)\phi_2(x)\phi_3(x)\d x&=&\frac{1}{(2\pi)^{2d}}\int_{\R^{2d}}
\widehat{\phi}_1(z)\widehat{\phi}_2(z_1)\overline{\widehat{\phi}_3(z+z_1)}\d
z\d z_1,\label{s2-13}
\eeqlb
(the inverse Fourier transform) and moreover, by the Riemann-Lebesgue Lemma,
$\widehat{\phi}_1(z)$ is bounded and goes to $0$ as $|z|\to\infty$.
It follows from (\ref{Eq:ChF}) and Fubini's theorem that for any $t>0$,
\beqlb\label{s2-14}
\widehat{T_t\phi}_1(z)=\widehat{\phi}_1(z)\e^{-t\sum_{k=1}^d|z_k|^{\alpha_k}}.
\eeqlb

From now on,  we define a sequence of random variables
$\tilde{X}_n$ in $\mathcal{S}'(\R^{d+1})$ as follows: For any
$n\geq0$ and $\psi\in\mathcal{S}(\R^{d+1})$,
 \beqlb\label{s2-15}
 \langle\tilde{X}_n, \psi\rangle=\int_0^1\langle X_n(t), \psi(\cdot, t)\rangle\d t.
 \eeqlb

\section{Laplace functionals of the occupation time fluctuations}

As in Bojdecki et al. \cite{BGT061, BGT062}, the method for proving
Theorems \ref{s2-thm-1}, \ref{s2-thm-3} and \ref{s2-thm-3} rely on
the Laplace functionals of the occupation time fluctuations. In this
section, we establish the main technical lemmas which will be needed
for the proofs in Section 4. For simplicity, we always assume
$\psi(x, t)= \phi(x)h(t)$, where $\phi\in\mathcal{S}(\R^{d})$ and
$h\in\mathcal{S}(\R)$ are nonnegative functions. Hence from
(\ref{s2-7}) and (\ref{s2-15}), we have that
 \beqlb\label{s3-1}
 \langle \tilde{X}_n, \psi\rangle=\frac{1}{F_n}\int_0^1\bigg(\int_0^{nt}
 \langle N_n(s), \phi\rangle\d s-\int_0^{nt} f_n(s)\langle\lambda, \phi\rangle\d
 s\bigg)h(t)\d t.\qquad
 \eeqlb
Let
 \beqlb\label{s3-2}
 \tilde{h}(s)=\int_s^1 h(t)\d t\;\;\text{and}\;\;\;\psi_n(x, s)
 =\frac{1}{F_n}\phi(x)\tilde{h}(\frac{s}{n}).
 \eeqlb
Then  (\ref{s3-1}) can be rewritten as
 \beqlb\label{s3-3}
 \langle \tilde{X}_n, \psi\rangle=\int_0^n\langle N_n(s),
 \psi_n(\cdot, s)\rangle\d s-\int_0^n f_n(s)\langle
 \lambda,\psi_n(\cdot, s)\rangle \d s.
 \eeqlb
Define
 \beqlb\label{s3-4}
 H_{n,\psi_n}(x, t, r)=\bfE_x\Bigg(\exp\bigg\{-\int_0^t\big\langle N_n(s),
 \psi_n(\cdot, r+s)\big\rangle\d s\bigg\}\Bigg).
 \eeqlb
Since $N_n(0)$ is a Poisson random measure with Lebesgue intensity
measure, it follows from (\ref{s3-3}) that
 \beqlb\label{s3-5}
 \bfE(\e^{-\langle \tilde{X}_n, \psi\rangle})
 =\exp\bigg\{\int_{\R^d}\big[H_{n,\psi_n}(x,n,0)-1\big]\d x
 +\int_{\R^d}\d x\int_0^n f_n(s)\psi_n(x, s)\d s\bigg\}.\qquad
 \eeqlb
Note that  for any $n\geq1$, $N_n$ is a Markov process. By using the
renewal argument, we can rewrite (\ref{s3-4}) as
 \beqlb
 \label{s3-1-1}
 H_{n,\psi_n}(x, t, r)
 &=&\e^{-\gamma t}\bfE_x\bigg\{\exp\bigg(-\int_0^t\psi_n(\vec{\xi}_n(s), r+s)\d
 s\bigg)\bigg\}\nonumber\\
 &&+\int_0^t\gamma\e^{-\gamma s}\bfE_x\bigg\{\exp\bigg(-\int_0^s\psi_n(\vec{\xi}_n(u), r+u)
 \d u\bigg)\nonumber  \\
 &&\times\bigg[\bfE_{\vec{\xi}_n(s)}\exp\bigg(-\int_0^{t-s}\big\langle N_n(u),
 \psi_n(\cdot, r+s+u)\big\rangle\d u\bigg)\bigg]
 ^{k_n(s)}\bigg\}\d s,\qquad
 \eeqlb
where $k_n(s)$ denotes the number of particles generated at the
first splitting time. Since the process $k_n$ is independent of the
process $\vec{\xi}_n$ and by the assumptions, for any $0\le z\le 1$
 $$
 \bfE(z^{k_n(s)})=g_n(z, s)=\Big(\frac{1+z^2}{2}\Big)\e^{-\delta_n s}+(1-\e^{-\delta_n s}),
 $$
(\ref{s3-1-1}) yields that
 \beqlb\label{s3-1-2}
 &&H_{n,\psi_n}(x, t, r)=\e^{-\gamma t}I_{n,\psi_n}(x, t, r)+
 \int_0^t\gamma\e^{-\gamma(t-s)}K_{n,\psi_n}(x, t-s, r, s)\d s,
 \eeqlb
where for any $x\in\R^d, t, r\geq0$,
 \beqnn
 I_{n,\psi_n}(x, t, r)&=&\bfE_x \exp\bigg(-\int_0^t\psi_n(\vec{\xi}_n(u), r+u)\d
 u\bigg),\label{s3-1-3}
  \\K_{n,\psi_n}(x, t, r, s)&=&\bfE_x \bigg[\exp\bigg(-\int_0^t
  \psi_n(\vec{\xi}_n(u), r+u)\d u\bigg)g_n\big(H_{n,\psi_n}(\vec{\xi}_n(t),s,r+t),
 t\big)\bigg].\label{s3-1-4}
 \eeqnn
Define $g'_{n,t}(z, t)=\frac{\partial g_n(z, t)}{\partial
t}=-\delta_n g_n(z, t)+\delta_n$ and $I_{n,\psi_n}(x, \d t,
r)=\frac{\partial I_{n,\psi_n}(x, t, r)}{\partial t}\d t$. By the
Feynman-Kac formula,
 \beqnn
 \frac{\partial I_{n,\psi_n}}{\partial t}&=&\Big(A+\frac{\partial}{\partial r}-\psi_n(x,
 r)\Big) I_{n,\psi_n},\label{s3-1-5}
 \\ \frac{\partial K_{n,\psi_n}}{\partial t}&=&\Big(A+\frac{\partial}{\partial r}-\psi_n(x,
 r)\Big)K_{n,\psi_n}
 \\&&+\bfE_x\Big(\e^{-\int_0^t\psi_n(\vec{\xi}_n(u),r+u)\d u}g'_{n,t}\big(H_{n,\psi_n}(\vec{\xi}_n(t),s,r+t),
 t\big)\Big)\nonumber
 \\&=&\Big(A+\frac{\partial}{\partial r}-\psi_n(x,
 r)\Big)K_{n,\psi_n}-\delta_n K_{n,\psi_n}+\delta_n I_{n,\psi_n}.\label{s3-1-6}
 \eeqnn
Therefore, (\ref{s3-1-2}) indicates that
 \beqnn
 \frac{\partial
H_{n,\psi_n}(x, t, r)}{\partial t}&=&\int_0^t\Big(A+\frac{\partial}{\partial r}-\psi_n(x,
 r)-\gamma\Big)\gamma\e^{-\gamma(t-s)}K_{n,\psi_n}(x, t-s, r, s)\d s\nonumber
 \\&+&\Big(A+\frac{\partial}{\partial r}-\psi_n(x,
 r)-\gamma\Big)\e^{-\gamma t}I_{n,\psi_n}(x, t, r)+\delta_n\int_0^t\gamma\e^{-\gamma s}I_{n,\psi_n}(x, s,
 r)\d s\nonumber
 \\&-&\delta_n\int_0^t\gamma\e^{-\gamma(t-s)}K_{n,\psi_n}(x, t-s,
 r,s)\d s+\gamma\bfE_x\Big[g_n\big(H_{n,\psi_n}(\vec{\xi}_n(0),t,r),
 0\big)\Big],\nonumber
 \eeqnn
which, combined with the fact that $g_n(z,
0)=\frac{1+z^2}{2}=:g(z)$, yields that
 \beqlb\label{s3-1-7}
 \frac{\partial
H_{n,\psi_n}(x, t, r)}{\partial
t}&=&\Big(A+\frac{\partial}{\partial r}-\psi_n(x,
 r)-\gamma-\delta_n\Big)H_{n,\psi_n}(x, t, r)+\gamma
 g(H_{n,\psi_n}(x,t,r))\nonumber
 \\&+&\delta_n\int_0^t \e^{-\gamma
 s}I_{n,\psi_n}(x, \d s, r)+\delta_n.
 \eeqlb
Let
 \beqlb\label{s3-1-8}
 V_{n,\psi_n}(x, t, r)=1-H_{n,\psi_n}(x, t, r).
 \eeqlb
Then $V_{n,\psi_n}(x, 0, r)=0$ and from (\ref{s3-1-7}), we have that
 \beqnn\label{s3-1-9}
\frac{\partial V_{n,\psi_n}(x, t, r)}{\partial t}&=&
 \Big(A+\frac{\partial}{\partial r}-\delta_n\Big)V_{n,\psi_n}(x, t, r)-\delta_n\int_0^t \e^{-\gamma
 s}I_{n,\psi_n}(x, \d s, r)
 \nonumber
 \\&&-\frac{\gamma}{2}V_{n,\psi_n}^2(x, t, r)+\psi_n(x, r)\big(1-V_{n,\psi_n}(x, t, r)\big),
 \eeqnn
where we have used the elementary fact that $g(1-z)-(1-z)=\frac{z^2}{2}$. Therefore,
 \beqlb\label{s3-1-10}
 &&V_{n,\psi_n}(x, t, r)=\int_0^t\e^{-\delta_ns}T_s\bigg\{\psi_n(\cdot, r+s)\big(1-V_{n,\psi_n}(\cdot, t-s, r+s)\big)
 -\frac{\gamma}{2}V_{n,\psi_n}^2(\cdot, t-s,
 r+s)\nonumber
 \\&&\qquad\qquad+\delta_n\int_0^{t-s}\e^{-\gamma u}\bfE_\cdot\Big[\e^{-\int_0^u\psi_n(\vec{\xi}_n(v), r+s+v)\d v}
 \psi_n(\vec{\xi}_n(u), r+s+u)\Big]
 \d u\bigg\}(x)\d s.\qquad
 \eeqlb
Some simple calculation shows
\beqlb\label{s3-1-12}
 &&\delta_n\int_0^t\e^{-\delta_ns}T_s\bigg(\int_0^{t-s}\e^{-\gamma u}T_u\psi_n(\cdot, r+s+u)\d u\bigg)(x)\d
 s\nonumber
\\&&\qquad\qquad\qquad=\int_0^t \big(f_n(s)-\e^{-\delta_n s}\big)T_{s}\psi_n(x, r+s)\d s.
  \eeqlb
In addition, by using (\ref{s2-1}), (\ref{s2-5}), (\ref{s3-4}), (\ref{s3-1-8}), and
the fact $1-\e^{-x}\leq x$ for all $x\in\R$, we derive that
\beqlb\label{s3-1-13}
V_{n,\psi_n}(x, t, r)\leq\int_0^t f_n(s)T_s\psi_n(\cdot, r+s)(x)\d s=:J_{n,\psi_n}(x, t,
 r).
 \eeqlb
It follows from (\ref{s3-1-10}), (\ref{s3-1-12}) and (\ref{s3-1-13}) that
 \beqlb\label{s3-1-14}
 J_{n,\psi_n}(x, t,
 r)- V_{n,\psi_n}(x, t,
 r)&=&\delta_n\int_0^t\e^{-\delta_n s}T_s\bigg(\int_0^{t-s}\e^{-\gamma u}
 \chi_{n,\psi_n}(\cdot, u, r+s)  \d u\bigg)(x)\d s\nonumber
 \\&+&\int_0^t\e^{-\delta_n s}T_s\bigg(\psi_n(\cdot, r+s)V_{n,\psi_n}(\cdot, t-s, r+s)\bigg)(x)\d
 s\nonumber
 \\&+&\int_0^t\e^{-\delta_n s}T_s\bigg(\frac{\gamma}{2}V_{n,\psi_n}^2(\cdot, t-s, r+s)\bigg)(x)\d
 s,
 \eeqlb
where for any $x\in\R^d, u\geq0$ and $s\geq0$
 \beqlb\label{s3-1-15}
\chi_{n,\psi_n}(x, u,
s)=\bfE_x\Big[\Big(1-\e^{-\int_0^u\psi_n(\vec{\xi}_n(v), s+v)\d
v}\Big)\psi_n(\vec{\xi}_n(u), s+u)\Big].
 \eeqlb
 By (\ref{s2-6}), (\ref{s3-5}) and (\ref{s3-1-14}), we obtain that
  \beqlb\label{s3-1-16}
   \bfE(\e^{-\langle \tilde{X}_n, \psi\rangle})&=&\exp\bigg(\int_{\R^d}\d x
   \int_0^n f_n(s)T_s\psi_n(x, s)\d s
   -\int_{\R^d}\big[1-H_{n,\psi_n}(x,n,0)\big]\d x\bigg)\nonumber
 \\&=&\exp\bigg(\int_{\R^d}\big[J_{n,\psi_n}(x, n,
 0)- V_{n,\psi_n}(x, n, 0)\big]\d x\bigg)\nonumber
 \\&=&\exp\Big(I_1(n,\psi_n)+I_2(n,\psi_n)+I_3(n,\psi_n)\Big),
   \eeqlb
 where
   \beqlb
 I_1(n,\psi_n)&=&\frac{\gamma}{2}\int_{\R^d}\d x\int_0^n\e^{-\delta_n s}V_{n,\psi_n}^2(x, n-s, s)\d
 s;\label{s3-1-17}
 \\ I_2(n,\psi_n)&=&\int_{\R^d}\d x\int_0^n\e^{-\delta_n s}\psi_n(x,
  s)V_{n,\psi_n}(x, n-s, s)\d s;\label{s3-1-18}
 \\ I_3(n,\psi_n)&=&\delta_n\int_{\R^d}\d x\int_0^n\e^{-\delta_n s}\d s\int_0^{n-s}
 \e^{-\gamma u}\chi_{n,\psi_n}(x, u, s) \d
  u.\label{s3-1-19}
  \eeqlb

Below, we evaluate the limits of $I_1(n,\psi_n)$, $I_2(n,\psi_n)$ and $I_3(n,\psi_n)$,
respectively. Recall that $H$ is the $d\times d$ diagonal matrix
$(1/\alpha_k)_{1\leq k\leq d}$. For all $t>0$ and $y=t^Hz$, $\d y=t^{\bar{\alpha}}\d z$.

We first study the limit of $I_1(n,\psi_n)$. By using (\ref{s3-1-17}), we can write
  \beqlb\label{s3-2-1}
  I_1(n,\psi_n)=I_{11}(n,\psi_n)+I_{12}(n,\psi_n),
  \eeqlb
where
 \beqlb
 I_{11}(n,\psi_n)&=&\frac{\gamma}{2}\int_{\R^d}\d x\int_0^n\e^{-\delta_n
 s}J_{n,\psi_n}^2(x, n-s, s)\d s,\label{s3-2-2}
 \\I_{12}(n,\psi_n)&=&\frac{\gamma}{2}\int_{\R^d}\d x\int_0^n\e^{-\delta_n
 s}(V_{n,\psi_n}^2(x, n-s, s)-J_{n,\psi_n}^2(x, n-s, s))\d
 s.\qquad\label{s3-2-3}
 \eeqlb

The following lemma determines the limit of $I_{11}(n,\psi_n)$.
\begin{lem}\label{s3-lem-1}
(1) If $\bar{\alpha} > 2$ and $F_n^2=n$, then as $n \to \infty$,
\beqlb\label{lem-1-1}
I_{11}(n,\psi_n)\to\frac{\gamma/2}{(2\pi)^{d}}\int_{\R^d}
\frac{|\widehat{\phi}(z)|^2}{(\sum_{k=1}^d|z_k|^{\alpha_k})^2}\d z
 \int_0^1h(u)\d u\int_0^1h(v)\d
 v\int_0^{u\wedge v}\e^{-\theta s}\d s.\qquad
 \eeqlb
(2) If $\bar{\alpha}=2$ and $F_n^2=n\ln n$, then
 \beqlb\label{lem-1-2}
I_1(n,\psi_n)\to C\int_0^1 h(r)\d r\int_0^1 h(t)\d t
 \int_0^{t\wedge r}\e^{-\theta s}\d s\int_{\R^d}\varphi\Big(\theta, r-s,
 t-s, \sum_{k=1}^d|y_k|^{\alpha_k}\Big)\d  y,\qquad
 \eeqlb
 where $C=\frac{\gamma/2}{(2\pi)^d}(\int_{\R^d}\phi(x)\d x)^2$ and
 $\varphi(x, u, v, z)$ is defined by (\ref{def1}).

\noindent (3) If $\bar{\alpha}\in(1, 2)$ and $F_n=n^{(3-\bar{\alpha})/2}$,
then
\beqlb\label{lem-1-3}
I_1(n,\psi_n)\to\frac{\gamma}{2\pi^{d}}\prod_{k=1}^d\frac{\Gamma(1/\alpha_k)}
{\alpha_k}\bigg(\int_{\R^d}\phi(x)\d x\bigg)^2\int_0^1\int_0^1 h(s)h(t)C(s, t)\d s\d t,\qquad
\eeqlb
where $C(s, t)$ is as (\ref{s2-thm-3-2}).
\end{lem}
{\bf Proof.} Since all
components of $\vec{\xi}$ are symmetric stable L\'{e}vy processes
and independent of each other, we have
 \beqlb\label{lem-1-4}
 \int_{\R^d} T_sf(x)g(x)\d x=\int_{\R^d}f(x)T_sg(x)\d x
 \eeqlb
for any $s\geq 0$ and bounded measurable functions $f$ and $g$.
By using (\ref{s2-4}), (\ref{s2-12}), (\ref{s2-14}), (\ref{s3-2}),
(\ref{s3-1-13}) and (\ref{s3-2-2}), we derive that
 \beqlb\label{lem-1-5}
 I_{11}(n,\psi_n)&=&\frac{\gamma}{2}\int_{\R^d}\d x\int_0^n\bigg(\e^{-\delta_n
 s}\int_0^{n-s}\e^{-\delta_n u}\bar{f}_n(u)T_u\psi_n(x, s+u)\d
 u\nonumber
 \\&&\qquad\qquad\quad\times\int_0^{n-s}\e^{-\delta_n v}\bar{f}_n(v)T_v\psi_n(x, s+v)\d v\bigg)\d
 s\nonumber
 \\&=&\frac{\gamma/2}{F_n^2(2\pi)^{d}}\int_0^n\bigg(\int_0^{n-s}\bar{f}_n(u)\tilde{h}\Big(\frac{s+u}{n}\Big)\d
 u\int_0^{n-s}\e^{-\delta_n (u+v)}\bar{f}_n(v)\tilde{h}\Big(\frac{s+v}{n}\Big)\d v\nonumber
 \\&&\qquad\qquad\quad\times\int_{\R^d}|\widehat{\phi}(z)|^2
 \e^{-(u+v)\sum_{k=1}^d|z_k|^{\alpha_k}}\d z\bigg)\e^{-\delta_n
 s}\d
 s,
 \eeqlb
Equation (\ref{lem-1-5}) and the fact that $\bar{f}_n$ converges uniformly to $1$ imply
that
 \beqlb\label{lem-1-6}
 \lim_{n\to\infty}I_{11}(n,\psi_n)&=&\lim_{n\to\infty}
 \bigg\{\frac{n^3\gamma/2 }{(2\pi)^{d}F_n^2}\int_0^1\bigg[
 \int_0^{1-s}\tilde{h}(s+u)\d
 u\int_0^{1-s}\e^{-n\delta_n (u+v)}\tilde{h}(s+v)\d v\nonumber
 \\&&\qquad\qquad\quad\times\int_{\R^d}|\widehat{\phi}(z)|^2\e^{-n(u+v)
 \sum_{k=1}^d|z_k|^{\alpha_k}}\d z\bigg]\e^{-n\delta_n s}\d
 s\bigg\}.
 \eeqlb

 {\bf (1)} Suppose $\bar{\alpha} >2$ and $F_n^2=n$. It follows from
(\ref{lem-1-6}) that
 \beqlb\label{lem-1-7}
 \lim_{n\to\infty}I_{11}(n,\psi_n)&=&\lim_{n\to\infty}\bigg\{\frac{\gamma/2}
 {(2\pi)^{d}}\int_{\R^d}|\widehat{\phi}(z)|^2\d z\int_0^1\bigg[\int_s^1
 n\e^{-n(u-s)(\delta_n+\sum_{k=1}^d|z_k|^{\alpha_k})}\tilde{h}(u)\d u\nonumber
 \\&&\qquad\times\int_s^1 n\e^{-n(v-s)(\delta_n+\sum_{k=1}^d|z_k|^{\alpha_k})}
 \tilde{h}(v)\d v\bigg]\e^{-n\delta_n s}\d
 s\bigg\}. \eeqlb
Since $\tilde{h}(t)$ is a continuous and bounded function, it is
easy to check that for every $z \in \R^d\backslash\{0\}$, as $n \to \infty$,
 $$\frac{\tilde{h}(s)}{\sum_{k=1}^d|z_k|^{\alpha_k}}\geq\int_s^1\tilde{h}(u)
 \e^{-n(u-s)(\delta_n+\sum_{k=1}^d|z_k|^{\alpha_k})}n\d
 u\to\frac{\tilde{h}(s)}{\sum_{k=1}^d|z_k|^{\alpha_k}}.$$
By using (\ref{lem-1-7}) and the dominated convergence theorem, we deduce that as
$n\to\infty$,
 \beqlb\label{lem-1-8}
 I_{11}(n,\psi_n)&\to&\frac{\gamma/2}{(2\pi)^d}\int_{\R^d}
 \Big(\frac{|\widehat{\phi}(z)|}{\sum_{k=1}^d|z_k|^{\alpha_k}}\Big)^2\d z\int_0^1\e^{-\theta s}
[\tilde{h}(s)]^2\d s.
 \eeqlb
Notice that, by Remark  \ref{s2-rm-1}, the right hand side of
(\ref{lem-1-8}) is finite. Substituting $\tilde{h}(s)=\int_s^1 h(t)\d t$
into (\ref{lem-1-8}) yields (\ref{lem-1-1}).

 {\bf (2)} Next we consider the case $\bar{\alpha}=2$ and $F_n^2=n\ln n$. Define
 \beqlb\label{lem-1-9}
 \bar{I}_{11}(n,\psi_n)&:=&\frac{n^3\gamma/2 }{(2\pi)^{d}F_n^2}\int_0^1\bigg[
 \int_0^{1-s}\tilde{h}(s+u)\d
 u\int_0^{1-s}\e^{-\theta (u+v)}\tilde{h}(s+v)\d v\nonumber
 \\&&\qquad\times\int_{\R^d}|\widehat{\phi}(z)|^2\e^{-n(u+v)\sum_{k=1}^d|z_k|^{\alpha_k}}\d z\bigg]\e^{-\theta s}\d
 s.
\eeqlb
Since $n\delta_n\to\theta$, for sufficiently large $n$ and all $t\in[0, 1]$,
 \beqnn
 \e^{-\theta t}(1-|\theta-n\delta_n|)\leq \e^{-\theta t}\e^{-|\theta-n\delta_n|t}
 &\leq&\e^{-\theta t}\e^{(\theta-n\delta_n)t}
 \\&=&\e^{-n\delta_n t}\leq\e^{-\theta t}\e^{|\theta-n\delta_n|t}\leq \e^{-\theta t}(1+3|\theta-n\delta_n|),
 \eeqnn
we see from  (\ref{lem-1-6}) that, in order to prove (\ref{lem-1-2}), it suffices to prove that
$\bar{I}_{11}(n,\psi_n)$ converges to
 \beqlb\label{lem-1-10}
 C\int_0^1 h(r)\d r\int_0^1 h(t)\d t  \int_0^{t\wedge r}\e^{-\theta s} \d s
 \int_{\R^d}\varphi\Big(\theta, r-s, t-s, \sum_{k=1}^d|y_k|^{\alpha_k}\Big)\d
 y,
 \eeqlb
 where $C=\frac{\gamma/2}{(2\pi)^d}(\int_{\R^d}\phi(x)\d x)^2$ and $\varphi$ is
 defined by (\ref{def1}). Substituting $\tilde{h}(u)=\int_u^1h(t)\d t$ and $F_n^2=n\ln n$
into (\ref{lem-1-9}) and changing the order of integration, we
obtain that
 \beqlb\label{lem-1-11}
 \lim_{n\to\infty}\bar{I}_{11}(n,\psi_n)&=&\lim_{n\to\infty}\bigg\{\frac{\gamma/2}{(2\pi)^{d}}
 \int_0^1 \e^{-\theta s}\d s\int_s^1\d u\int_u^1 h(t)\d t\int_s^1\d v\int_v^1h(r)\d r\nonumber
 \\&&\qquad\quad\times\int_{\R^d}|\widehat{\phi}(z)|^2\e^{-(u+v-2s)(\theta+n\sum_{k=1}^d|z_k|^{\alpha_k})}\d z\bigg\}\nonumber
 \\&=&\lim_{n\to\infty}\bigg\{\frac{\gamma/2}{(2\pi)^{d}}
 \int_0^1 \e^{-\theta s}\d s\int_s^1 h(t)\d t\int_s^1h(r)\d r\int_{\R^d}|\widehat{\phi}(z)|^2\d z\nonumber
 \\&&\qquad\quad\times\int_s^t\int_s^r\e^{-(u+v-2s)(\theta+n\sum_{k=1}^d|z_k|^{\alpha_k})}\d u\d v\bigg\}\nonumber
  \\&=&\lim_{n\to\infty}\frac{\gamma/2}{(2\pi)^{d}}\int_0^1 h(r)\d r\int_0^1
 h(t)\d t\int_0^{r\wedge t}\e^{-\theta s}W_{n, r, t, s}(\theta)\d
 s,\qquad
 \eeqlb
where for any $x\geq0$ and $0 \le s \le r \wedge t$,
 \beqlb\label{lem-1-12}
 &&W_{n, r, t,s}(x)=\int_{\R^d}|\widehat{\phi}(z)|^2
 \frac{\varpi(x,r-s,t-s,\sum_{k=1}^d|z_k|^{\alpha_k},n)}{\ln n}\d z,\qquad\quad
 \eeqlb
and for $x, u, v\geq 0$ and $y, w>0$,
 $$\varpi(x,u,v,w,y)=y^2\frac{(1-\e^{-u(x+yw)})(1-\e^{-v(x+yw)})}{(x+yw)^2}.$$
Since the function $r \mapsto (1-\e^{-r})/r$  is decreasing in
$r\in(0,+\infty)$, we have that
 \beqlb\label{lem-1-13}
 \varpi(x,u,v,w,y)\leq \varpi(0,1,1,w,y)
 \eeqlb for all $u, v\in[0,
1]$. By applying  L'H\^{o}pital's rule and the dominated convergence
theorem, we derive from (\ref{lem-1-12}) that \beqnn
\lim_{n\to\infty}W_{n, r, t, s}(\theta)
=\lim_{n\to\infty}\int_{\R^d}|\widehat{\phi}(z)|^2\bar{\varphi}_n(\theta,
r, t, s, z)\d z,\qquad \eeqnn where \beqnn \bar{\varphi}_n(\theta,
r, t, s,
z)&=&n\varpi_y'\Big(\theta,r-s,t-s,\sum_{k=1}^d|z_k|^{\alpha_k},n\Big)
=n^2\varphi\Big(\theta, r-s,t-s,n\sum_{k=1}^d|z_k|^{\alpha_k}\Big),
\eeqnn and $\varphi$ is defined by (\ref{def1}). By using the
substitution $y=n^{H}z$, we get that
 $$
 \int_{\R^d}|\widehat{\phi}(z)|^2\bar{\varphi}_n(\theta, r, t, s, z)\d z
 =\int_{\R^d}|\widehat{\phi}(n^{-H}y)|^2\varphi\Big(\theta,r-s, t-s, \sum_{k=1}^d|y_k|^{\alpha_k}\Big)\d y.$$
Therefore,
  \beqlb\label{lem-1-14}
 \lim_{n\to\infty}W_{n, r, t, s}(\theta)
 =|\widehat{\phi}(0)|^2\int_{\R^d}\varphi\Big(\theta,r-s, t-s, \sum_{k=1}^d|y_k|^{\alpha_k}\Big)\d y.\qquad
 \eeqlb
On the other hand, let
 \beqnn
 \tilde{W}(n)=W_{n, 1, 1, 0}(0)=\int_{\R^d}|\widehat{\phi}(z)|^2
 \frac{\varpi(0,1,1,\sum_{k=1}^d|z_k|^{\alpha_k},n)}{\ln n}\d z.
 \eeqnn
Then (\ref{lem-1-14}) yields that
 \beqlb\label{lem-1-15}
 \lim_{n\to\infty}\tilde{W}(n)&=&|\widehat{\phi}(0)|^2\int_{\R^d}
 \varphi\Big(0,1, 1, \sum_{k=1}^d|y_k|^{\alpha_k}\Big)\d y\nonumber
 \\&=&2|\widehat{\phi}(0)|^2\int_{\R^d}\frac{\e^{-\sum_{k=1}^d|z_k|^{\alpha_k}}
 (1-\e^{-\sum_{k=1}^d|z_k|^{\alpha_k}})}{\sum_{k=1}^d|z_k|^{\alpha_k}}\d z<\infty.
 \eeqlb
Hence $\{\tilde{W}(n)\}_{n\geq 1}$ is bounded. Combining
(\ref{lem-1-12}), (\ref{lem-1-13}) and (\ref{lem-1-15}) gives
 $$
 |W_{n, r, t, s}(\theta)|\leq |\tilde{W}(n)|\leq K,
 $$
for some constant $K<\infty$. By the dominated convergence theorem, we derive (\ref{lem-1-10})
from (\ref{lem-1-11})  and (\ref{lem-1-14}).

{\bf (3)} To consider the case $\bar{\alpha} \in (1, 2)$, we
substitute $y=(nu+nv)^{H}z$ into (\ref{lem-1-6})  to obtain
that
 \beqlb\label{lem-1-16}
\lim_{n\to\infty}I_{11}(n,\psi_n)&=&\lim_{n\to\infty} \bigg\{\frac{
n^{3-\bar{\alpha}}\gamma/2}{(2\pi)^{d}F_n^2}\int_{\R^d}|\widehat{\phi}((nu+nv)^{-H}y)|^2
\e^{-\sum_{k=1}^d|y_k|^{\alpha_k}}\d y\int_0^1\bigg[\int_0^{1-s}\d
 u\nonumber
 \\&&\times\tilde{h}(s+u)\int_0^{1-s}\e^{-n\delta_n (u+v)}\tilde{h}(s+v)
 \frac{\d v}{(u+v)^{\bar{\alpha}}}\bigg]\e^{-n\delta_n s}\d
 s\bigg\}.
 \eeqlb
Since $(nu+nv)^{-H}y\to 0$ for every $y$, substituting
$F_n=n^{(3-\bar{\alpha})/2}$ into (\ref{lem-1-16}), we see that
$I_{11}(n,\psi_n)$ converges to
 \beqlb\label{lem-1-17}
&&\frac{\gamma/2
}{(2\pi)^{d}}\int_{\R^d}|\widehat{\phi}(0)|^2\e^{-\sum_{k=1}^d|y_k|^{\alpha_k}}\d
y\int_0^1\e^{\theta s}\bigg[\int_s^1\e^{-\theta u}\tilde{h}(u)\d
 u\int_s^{1}\e^{-\theta v}\frac{\tilde{h}(v)\d
v}{(u+v-2s)^{\bar{\alpha}}}\bigg]\d s\nonumber
\\&&=\frac{\gamma/2}{\pi^{d}}|\widehat{\phi}(0)|^2\prod_{k=1}^d
\frac{\Gamma(1/\alpha_k)}{\alpha_k}\int_0^1\e^{-\theta u}\tilde{h}(u)\d u
\int_0^1\e^{-\theta v}\tilde{h}(v)\d v\int_0^{v\wedge u}\frac{\e^{\theta s}\d s}
{(u+v-2s)^{\bar{\alpha}}}.
\eeqlb
Substituting $\tilde{h}(t)=\int_t^1h(s)\d s$ into (\ref{lem-1-17}) yields (\ref{lem-1-3}).\qed

\medskip

For $I_{12}(n,\psi_n)$ in (\ref{s3-2-3}) we have

\begin{lem}\label{s3-lem-2}
If $\bar{\alpha}>1$ and $F_n$ takes values according to Theorem \ref{s2-thm-1},
\ref{s2-thm-2} and \ref{s2-thm-3}, respectively, then
\beqlb\label{lem-2-1}
\lim_{n\to\infty}I_{12}(n,\psi_n)=0.
\eeqlb
\end{lem}
{\bf Proof}. By (\ref{s3-1-13}) and (\ref{s3-1-14}), we have that
 \beqnn
 &&J^2_{n,\psi_n}(x, n-s, s)-V^2_{n,\psi_n}(x, n-s, s)
 \\&&\qquad\leq 2J_{n,\psi_n}(x, n-s, s)\bigg[\delta_n\int_0^{n-s}\e^{-\gamma u}T_u\psi_n(x,
 s+u)\d u\int_0^u\e^{-(\delta_n-\gamma)v }\d v\nonumber
 \\&&\qquad+\int_0^{n-s}\e^{-\delta_n u}T_u\big(\psi_n(\cdot, s+u)J_{n,\psi_n}(\cdot, n-s-u,
 s+u)\big)(x)\d u
 \\&&\qquad+\frac{\gamma}{2}\int_0^{n-s}\e^{-\delta_n u}T_uJ_{n,\psi_n}^2(x, n-s-u, s+u)\d
 u\bigg].
 \eeqnn
This and (\ref{s3-2-3}) imply that
 \beqlb\label{lem-2-2}
 |I_{12}(n,\psi_n)|\leq\gamma I_{121}(n,\psi_n)+\gamma
 I_{122}(n,\psi_n)+\frac{\gamma^2}{2}I_{123}(n,\psi_n),
 \eeqlb
where
 \beqlb
 I_{121}(n,\psi_n)&=&\int_{\R^d}\d x\int_0^n\e^{-\delta_n
 s}J_{n,\psi_n}(x, n-s, s)\d
 s\nonumber
 \\&&\times\int_0^{n-s}\e^{-\delta_n u}(\bar{f}_n(u)-1)T_u\psi_n(x,
 s+u)\d u,\qquad\label{lem-2-3}
 \\I_{122}(n,\psi_n)&=&\int_{\R^d}\d x\int_0^n\e^{-\delta_n
 s}J_{n,\psi_n}(x, n-s, s)\d
 s\nonumber
 \\&&\times\int_0^{n-s}\e^{-\delta_n u}T_u\big(\psi_n(\cdot, s+u)J_{n,\psi_n}(\cdot, n-s-u,
 s+u)\big)(x)\d u  \qquad\label{lem-2-4}
 \eeqlb
and
 \beqlb
 I_{123}(n,\psi_n)&=&\int_{\R^d}\d x\int_0^n\e^{-\delta_n
 s}J_{n,\psi_n}(x, n-s, s)\d
 s\nonumber
 \\&&\times\int_0^{n-s}\e^{-\delta_n u}T_uJ_{n,\psi_n}^2(x, n-s-u, s+u)\d
 u.\label{lem-2-5}
 \eeqlb

Substituting  (\ref{s3-1-13}) and (\ref{s3-2}) into (\ref{lem-2-3})
gives that
 \beqlb\label{lem-2-6}
I_{121}(n,\psi_n)&=&\frac{1}{F_n^2}\int_0^n\e^{-\delta_n s}\d
s\int_0^{n-s}(\bar{f}_n(u)-1)\e^{-\delta_n
u}\tilde{h}_n\Big(\frac{s+u}{n}\Big)\d u\nonumber
\\
&&\times\int_0^{n-s}\e^{-\delta_n v}\bar{f}_n(v)\tilde{h}_n\Big(\frac{s+v}{n}\Big)\d v
\int_{\R^d}T_u\phi(x)T_v\phi(x)\d x.
 \eeqlb
Note that  from (\ref{s2-12}), (\ref{s2-14}) and (\ref{lem-1-4}) we have that
$$
\int_{\R^d}T_u\phi(x)T_v\phi(x)\d x=\frac{1}{(2\pi)^d}
 \int_{\R^d}|\widehat{\phi}(z)|^2\e^{-(u+v)\sum_{k=1}^d|z_k|^{\alpha_k}}\d z.
$$
Comparing (\ref{lem-2-6}) with (\ref{lem-1-5}), we find that as
$n\to\infty$,
\beqlb\label{lem-2-7}
I_{121}(n,\psi_n)\leq
\frac{2\delta_n}{\gamma(\gamma-\delta_n)}I_{11}(n,\psi_n),
\eeqlb
where we have used the fact
$0<(\bar{f}_n(u)-1)/\bar{f}_n(u)\leq\delta_n/(\gamma-\delta_n)$ which follows from (\ref{s2-4}).

In addition, substituting (\ref{s2-4}) and (\ref{s3-1-13}) into (\ref{lem-2-4}) we
obtain that
 \beqlb\label{lem-2-8}
 I_{122}(n,\psi_n)&=&\int_0^n\e^{-\delta_n s}\d s\int_0^{n-s}\e^{-\delta_n
 u}\d u\int_0^{n-s}\bar{f}_n(v)\e^{-\delta_n v}T_v\psi_n(x, s+v)\d
 v\nonumber
 \\&&\times\int_{\R^d}T_u\Big(\psi_n(\cdot, s+u)\int_0^{n-s-u}\bar{f}_n(w)
 \e^{-\delta_n w}T_w\psi_n(\cdot, s+u+w)\d
 w\Big)(x)\d x\nonumber
 \\&=&\int_0^n\e^{-\delta_n s}\d s\int_0^{n-s}\e^{-\delta_n
 u}\d u\int_0^{n-s}\bar{f}_n(v)\e^{-\delta_n v}\d
 v\int_0^{n-s-u}\bar{f}_n(w)\e^{-\delta_n w}\d w\nonumber
 \\&&\times\int_{\R^d}T_v\psi_n(x, s+v)T_u\big(\psi_n(\cdot, s+u)T_w\psi_n(\cdot, s+u+w)\big)(x)\d
 x.
 \eeqlb
Using (\ref{s2-13}), (\ref{s2-14}), (\ref{s3-2}) and (\ref{lem-1-4}),
we derive from (\ref{lem-2-8}) that
 \beqnn
 I_{122}(n,\psi_n)&&=\int_0^n\e^{-\delta_n s}\d s\int_0^{n-s}\e^{-\delta_n
 u}\d u\int_0^{n-s}\bar{f}_n(v)\e^{-\delta_n v}\d
 v\int_0^{n-s-u}\bar{f}_n(w)\e^{-\delta_n w}\d w\nonumber
 \\&&\qquad\qquad\times\int_{\R^d}\psi_n(x, s+u)
 T_w\psi_n(x, s+u+w)T_{v+u}\psi_n(x, s+v)\d  x\nonumber \\
&&=\frac{1}{(2\pi)^{2d}F_n^3}\int_{\R^{2d}}\widehat{\phi}(z)\widehat{\phi}(z')
\overline{\widehat{\phi}(z+z')}\d z\d z'\int_0^n\d s
\int_0^{n-s}\e^{-\delta_n (s+u)}\tilde{h}\Big(\frac{s+u}{n}\Big)\d u
\nonumber
\\&&\qquad\qquad\times \int_0^{n-s}\bar{f}_n(v)\e^{-\delta_n v}\tilde{h}\Big(\frac{s+v}{n}\Big)\d
v\int_0^{n-s-u}\bigg[\bar{f}_n(w)\e^{-\delta_n w}\tilde{h}\Big(\frac{s+u+w}{n}\Big)\nonumber
\\&&\qquad\qquad\times \e^{-(u+v)\sum_{k=1}^d|z_k|^{\alpha_k}}
\e^{-w\sum_{k=1}^d|z'_k|^{\alpha_k}}\bigg]\d w.
\eeqnn
Note that $\{\bar{f}_n\}_n$ and $|\widehat{\phi}|$ are bounded
and, for  $h\in\mathcal{S}(\R)$, $\tilde{h}$ is bounded as well.
There exists a constant $K>0$ such that
\beqlb\label{lem-2-9}
 I_{122}(n,\psi_n)&\leq&\frac{Kn^4}{F_n^3}\int_{\R^{2d}}|\widehat{\phi}(z)\widehat{\phi}(z')|\d z\d z'
 \int_0^1\e^{-n\delta_n s}\d s
\int_0^{1-s}\d u\int_0^{1-s}\d
 v \nonumber
\\&&\times \int_0^{1-s-u}\e^{-n(u+v)\sum_{k=1}^d|z_k|^{\alpha_k}}\e^{-nw\sum_{k=1}^d|z'_k|^{\alpha_k}}\d w.
\eeqlb

Furthermore, substituting (\ref{s2-4}) and (\ref{s3-1-13}) into (\ref{lem-2-5}), we
obtain that
 \beqnn
I_{123}(n,\psi_n)&=&\int_0^n\e^{-\delta_n s}\d
s\int_0^{n-s}\e^{-\delta_n u}\d u\int_0^{n-s-u}\e^{-\delta_n
t}\bar{f}_n(t)\d t\nonumber
\\&&\times\int_0^{n-s-u}\e^{-\delta_n t'}\bar{f}_n(t')\d
t'\int_0^{n-s}\e^{-\delta_n v}\bar{f}_n(v)\d
v
\\&&\times\int_{\R^d}T_t\psi_n(x, s+u+t)T_{t'}\psi_n(x,
s+u+t')T_{u+v}\psi_n(x, s+u)\d x.\nonumber
 \eeqnn
By the same argument used to get (\ref{lem-2-9}),  we can find $K>0$ such that
 \beqlb\label{lem-2-10}
I_{123}(n,\psi_n)&\leq&\frac{Kn^5}{F_n^3}\int_{\R^{2d}}|\widehat{\phi}(z+z')
\widehat{\phi}(z')\widehat{\phi}(z)|\d z\d z'
\int_0^1\e^{-n\delta_n s}\d s\int_0^{1-s}\e^{-nu\sum_{k=1}^d|z_k|^{\alpha_k}}\d
u\nonumber
\\&&\times\int_0^{1-s}\e^{-nv\sum_{k=1}^d|z_k|^{\alpha_k}}\d v\int_0^{1-s-u}
\e^{-nt'\sum_{k=1}^d|z'_k|^{\alpha_k}}\d
t'\nonumber
\\&&\times\int_0^{1-s-u}\e^{-nt\sum_{k=1}^d|z_k+z'_k|^{\alpha_k}}\d t.
\eeqlb

Note that Lemma \ref{s3-lem-1} and (\ref{lem-2-7}) imply
\beqnn
\lim_{n\to\infty}I_{121}(n,\psi_n)=0.
\eeqnn
Thus  we see from (\ref{lem-2-2}) that, in order to prove (\ref{lem-2-1}), it suffices to
show that $I_{122}(n,\psi_n)$ and $I_{123}(n,\psi_n)$ all converge to $0$ as $n\to\infty$.
Below we divide the proof of these facts into three cases.

\medskip

{\bf Case (1)} $\bar{\alpha}\in(2, \infty)$ and $F_n^2=n$. From
(\ref{lem-2-9}), we have that for some $K>0$
\beqnn
I_{122}(n,\psi_n)&\leq&\frac{K}{\sqrt{n}}\int_{\R^{2d}}\frac{|\widehat{\phi}(z)|}
{\big(\sum_{k=1}^d|z_k|^{\alpha_k}\big)^2}
\frac{|\widehat{\phi}(z')|}{\sum_{k=1}^d|z'_k|^{\alpha_k}}\d z\d z'.
\eeqnn
By Remark \ref{s2-rm-1} the last integral is finite. Hence we have
\beqlb\label{lem-2-12}
\lim\limits_{n\to\infty}I_{122}(n,\psi_n)=0,
\eeqlb

From (\ref{lem-2-10}), it follows that
\beqlb\label{lem-2-17}
I_{123}(n,\psi_n)&\leq&\frac{K}{\sqrt{n}}\int_{\R^{2d}}\frac{|\widehat{\phi}(z+z')|}
{\sum_{k=1}^d|z_k+z'_k|^{\alpha_k}}
 \frac{|\widehat{\phi}(z')|}{\sum_{k=1}^d|z'_k|^{\alpha_k}}
 \frac{|\widehat{\phi}(z)|}{(\sum_{k=1}^d|z_k|^{\alpha_k})^2}\d z\d
 z'.\qquad
 \eeqlb
By the same argument as those used in \cite[p.27]{BGT062}, we can verify
that
 $$\int_{\R^{2d}}\frac{|\widehat{\phi}(z+z')|}{\sum_{k=1}^d|z_k+z'_k|^{\alpha_k}}
 \frac{|\widehat{\phi}(z')|}{\sum_{k=1}^d|z'_k|^{\alpha_k}}
 \frac{|\widehat{\phi}(z)|}{(\sum_{k=1}^d|z_k|^{\alpha_k})^2}\d z\d
 z'<\infty.$$
Hence (\ref{lem-2-17}) implies that
 \beqlb\label{lem-2-15}
\lim_{n \to \infty} I_{123}(n,\psi_n) = 0.
 \eeqlb

{\bf Cases (2)} $\bar{\alpha}=2$ and $F_n^2=n\ln n$.  (\ref{lem-2-9}) implies that
for some $K>0$,
 \beqlb\label{lem-2-16}
 I_{122}(n,\psi_n)\leq\frac{K}{(n\ln n)^{1/2}}\int_{\R^{d}}
 \frac{|\widehat{\phi}(z)|(1-\e^{-n\sum_{k=1}^d|z_k|^{\alpha_k}})^2}
 {\big(\sum_{k=1}^d|z_k|^{\alpha_k}\big)^2\ln n}\d z
 \int_{\R^d}\frac{|\widehat{\phi}(z')|\d
 z'}{\sum_{k=1}^d|z'_k|^{\alpha_k}}.\qquad
 \eeqlb
By applying L'H\^{o}pital's rule and substituting $y=n^Hz$, we derive that as $n \to \infty$,
$$
\int_{\R^{d}}  \frac{|\widehat{\phi}(z)|(1-\e^{-n\sum_{k=1}^d|z_k|^{\alpha_k}})^2}
 {\big(\sum_{k=1}^d|z_k|^{\alpha_k}\big)^2\ln n}\d z \to
 2|\widehat{\phi}(0)|\int_{\R^{d}}\frac{\e^{-\sum_{k=1}^d|y_k|^{\alpha_k}}
 (1-\e^{-\sum_{k=1}^d|y_k|^{\alpha_k}})}{\sum_{k=1}^d|y_k|^{\alpha_k}}\d y.
 $$
By Remark \ref{s2-rm-1}, the last integral and $\int_{\R^d}\frac{|\widehat{\phi}(z')|}
{\sum_{k=1}^d|z'_k|^{\alpha_k}}\d z'$  are finite. Therefore (\ref{lem-2-16})
implies (\ref{lem-2-12}).

Substituting $F_n^2=n\ln n$ into (\ref{lem-2-10}), we get that for some constant $K>0$,
 \beqnn
  I_{123}(n,\psi_n)&\leq&\frac{K}{ n^{1/2}(\ln n)^{3/2}}
  \int_{\R^{2d}}\bigg[\frac{(1-\e^{-n\sum_{k=1}^d|z_k+z'_k|^{\alpha_k}})
  (1-\e^{-n\sum_{k=1}^d|z_k|^{\alpha_k}})}
 {\sum_{k=1}^d|z_k+z'_k|^{\alpha_k}\sum_{k=1}^d|z'_k|^{\alpha_k}}
 \nonumber
 \\&&\qquad\qquad\qquad\quad\times\frac{(1-\sum_{k=1}^d|z_k|^{\alpha_k})^2}
 {(\sum_{k=1}^d|z_k|^{\alpha_k})^2}
 |\widehat{\phi}(z)||\widehat{\phi}(z')||\widehat{\phi}(z+z')|\bigg]\d z\d
 z'.
 \eeqnn
 Furthermore, by using the  inequality $1-\e^{-x}\leq
 x^{1/8}$ for $x\geq 0$ and the similar argument to that in \cite[p.29 lines 8-15]{BGT062},
  we arrive at (\ref{lem-2-15}). The details are omitted.

\medskip

{\bf Case (3)} $\bar{\alpha}\in(1, 2)$ and $F_n^2=n^{3-\bar{\alpha}}$. It follows from
the boundedness of $|\widehat{\phi}(z)|$  and (\ref{lem-2-9}) that
\beqlb\label{lem-2-11}
I_{122}(n,\psi_n)\leq\frac{Kn}{F_n^3}\int_{\R^{2d}}\frac{|\widehat{\phi}(z')|}
{\sum_{k=1}^d|z'_k|^{\alpha_k}} \bigg[\frac{1-\e^{-n\sum_{k=1}^d|z_k|^{\alpha_k}}}
{\sum_{k=1}^d|z_k|^{\alpha_k}}\bigg]^2\d
 z\d z'.\qquad
 \eeqlb
Substituting $y=n^H z$ in (\ref{lem-2-11}) yields
\beqlb\label{lem-2-11b}
I_{122}(n,\psi_n)&\leq&\frac{K}{F_n}\int_{\R^{2d}}\frac{|\widehat{\phi}(z')|\d
z'}{\sum_{k=1}^d|z'_k|^{\alpha_k}}
\bigg[\frac{1-\e^{-\sum_{k=1}^d|y_k|^{\alpha_k}s}}{\sum_{k=1}^d|y_k|^{\alpha_k}}\bigg]^2\d
y.
\eeqlb
Note that, for $\bar{\alpha}\in(1, 2)$, Lemma \ref{s2-lm-1} and Remark
\ref{s2-rm-1} imply
$$
\int_{\R^d}\frac{|\widehat{\phi}(z)|}{\sum_{k=1}^d|z_k|^{\alpha_k}}\d z\int_{\R^d}
 \bigg[\frac{1-\e^{-\sum_{k=1}^d|y_k|^{\alpha_k}}}{\sum_{k=1}^d|y_k|^{\alpha_k}}\bigg]^2\d
 y<\infty.$$
Thus, (\ref{lem-2-11b}) implies (\ref{lem-2-12}).

On the other hand, (\ref{lem-2-10}) implies that
\beqlb \label{lem-2-13}
I_{123}(n,\psi_n)&\leq&\frac{Kn}{F_n^3}\int_{\R^{2d}}\frac{|\widehat{\phi}(z+z')|}
{\sum_{k=1}^d|z_k+z'_k|^{\alpha_k}}\frac{|\widehat{\phi}(z')|}{\sum_{k=1}^d|z'_k|^{\alpha_k}}
|\widehat{\phi}(z)|\d z\d z'\nonumber
\\&& \qquad \qquad \times\int_0^1\bigg[\frac{1-\e^{-n(1-s)\sum_{k=1}^d|z_k|^{\alpha_k}}}
{\sum_{k=1}^d|z_k|^{\alpha_k}}\bigg]^2\d s \\
&\leq&\frac{Kn}{F_n^3}\int_{\R^{2d}}\frac{1}{\sum_{k=1}^d|z_k+z'_k|^{\alpha_k}}
\frac{1}{\sum_{k=1}^d|z'_k|^{\alpha_k}}
\bigg[\frac{1-\e^{-n\sum_{k=1}^d|z_k|^{\alpha_k}}}{\sum_{k=1}^d|z_k|^{\alpha_k}}\bigg]^2\d
z\d z'. \nonumber
 \eeqlb
Letting $y=n^{H}z$ and $y'=n^{H}z'$ and substituting
$F_n^2=n^{3-\bar{\alpha}}$ into (\ref{lem-2-13}) lead to
 \beqlb\label{lem-2-14}
I_{123}(n,\psi_n)&\leq&\frac{Kn^{1/2}}{n^{\bar{\alpha}/2}}\int_{\R^{2d}}
\frac{1}{\sum_{k=1}^d|y_k+y'_k|^{\alpha_k}}\frac{1}{\sum_{k=1}^d|y'_k|^{\alpha_k}}
\bigg[\frac{1-\e^{-\sum_{k=1}^d|y_k|^{\alpha_k}}}{\sum_{k=1}^d|y_k|^{\alpha_k}}\bigg]^2\d
y\d y'.\qquad
 \eeqlb
Since $\bar{\alpha}\in(1, 2)$, we can use the same argument as in
\cite[p.17]{BGT061} to verify that
 \beqnn
 \int_{\R^{2d}}\frac{1}{\sum_{k=1}^d|y_k+y'_k|^{\alpha_k}}\frac{1}{\sum_{k=1}^d|y'_k|^{\alpha_k}}
\bigg[\frac{1-\e^{-\sum_{k=1}^d|y_k|^{\alpha_k}}}{\sum_{k=1}^d|y_k|^{\alpha_k}}\bigg]^2\d
y\d y'<\infty.
 \eeqnn
Therefore, (\ref{lem-2-15}) follows from (\ref{lem-2-14}).
\qed

\medskip

For $I_2(n,\psi_n)$ in (\ref{s3-1-18}), we have the following lemma.

\begin{lem}\label{s3-lem-3}
(1) If $\bar{\alpha}>2$ and $F_n^2=n$, then
\beqlb\label{lem-3-1}
\lim_{n\to\infty}I_2(n,\psi_n)=\frac{1}{(2\pi)^{d}}\int_0^1h(t)\d t
\int_0^1 h(t')\d t'\int_0^{t\wedge t'}\e^{-\theta s}\d s
\int_{\R^d}\frac{|\widehat{\phi}(z)|^2}{\sum_{k=1}^d|z_k|^{\alpha_k}}\d z.
 \eeqlb
(2) If $\bar{\alpha}\in(1, 2]$ and $F_n$ takes values according to Theorems
\ref{s2-thm-2} and \ref{s2-thm-3}, respectively, then
 \beqnn
 \lim_{n\to\infty}I_{2}(n,\psi_n)=0.
 \eeqnn
\end{lem}
{\bf Proof.} To prove (\ref{lem-3-1}), we write
 \beqlb\label{lem-3-4}
 I_2(n,\psi_n)&=&I_{21}(n,\psi_n)-I_{22}(n,\psi_n),
 \eeqlb
where
 \beqnn
 I_{21}(n,\psi_n)&=&\int_{\R^d}\d x\int_0^n\e^{-\delta_n s}\psi_n(x,
 s)J_{n,\psi_n}(x, n-s, s)\d s,
 \\I_{22}(n,\psi_n)&=&\int_{\R^d}\d x\int_0^n\e^{-\delta_n s}\psi_n(x,
 s)\big(J_{n,\psi_n}(x, n-s, s)-V_{n,\psi_n}(x, n-s, s)\big)\d
 s.
 \eeqnn
By (\ref{s2-12}), (\ref{s2-14}), (\ref{s3-2}) and (\ref{s3-1-13}), $I_{21}(n,\psi_n)$
equals
\beqnn
 \frac{n^2}{F_n^2(2\pi)^{d}}\int_0^1\e^{-n\delta_n s}\d
 s \int_0^{1-s}f_n(nu)
 \tilde{h}(s)\tilde{h}(s+u)\d u\int_{\R^d}|\widehat{\phi}(z)|^2
 \e^{-nu\sum_{k=1}^d|z_k|^{\alpha_k}}\d z.
\eeqnn Substituting $\tilde{h}(t)=\int_t^1 h(s)\d s$ and
(\ref{s2-4}) into the above formula gives that
 \beqlb\label{lem-3-5}
 I_{21}(n,\psi_n)&=&\frac{n^2}{F_n^2(2\pi)^{d}}\int_0^1\e^{-n\delta_n s}\d
 s \int_0^{1-s}\e^{-n\delta_n u}\bar{f}_n(nu)\d u
 \int_s^1h(t)\d t\nonumber
 \\&&\qquad\quad \times\int_{s+u}^1 h(t')\d t'\int_{\R^d}|\widehat{\phi}(z)|^2
 \e^{-nu\sum_{k=1}^d|z_k|^{\alpha_k}}\d
 z\nonumber
 \\&=&\frac{n^2}{F_n^2(2\pi)^{d}}\int_0^1\e^{-n\delta_n s}\d
 s \int_s^1h(t)\d t\int_s^1 h(t')\d t'
 \nonumber
 \\&&\qquad\times\int_{\R^d}|\widehat{\phi}(z)|^2\d
 z\int_0^{t'-s}\bar{f}_n(nu)\e^{-nu(\delta_n+\sum_{k=1}^d|z_k|^{\alpha_k})}\d
 u.
 \eeqlb
Since $\bar{f}_n(u)\to 1$ uniformly as $n\to\infty$, (\ref{lem-3-5})
implies that
 \beqlb\label{lem-3-6}
 \lim_{n\to\infty}I_{21}(n,\psi_n)=\lim_{n\to\infty}\tilde{I}_{21}(n,\psi_n),
 \eeqlb
where
 \beqnn
 \tilde{I}_{21}(n,\psi_n) &=&\frac{n^2}{F_n^2(2\pi)^{d}}\int_0^1\e^{-n\delta_n
s}\d  s \int_s^1h(t)\d t\int_s^1 h(t')\d t'\nonumber
 \\&&\qquad\qquad\times \int_{\R^d}|\widehat{\phi}(z)|^2\d
 z\int_0^{t'-s}\e^{-nu(\delta_n+\sum_{k=1}^d|z_k|^{\alpha_k})}\d
 u
 \\&&=\frac{n}{F_n^2(2\pi)^{d}}\int_0^1\e^{-n\delta_n
s}\d  s \int_s^1h(t)\d t\int_s^1 h(t')\d t'\nonumber
 \\&&\qquad\qquad\times\int_{\R^d}|\widehat{\phi}(z)|^2\frac{1-\e^{-n(t'-s)(\delta_n+
 \sum_{k=1}^d|z_k|^{\alpha_k})}}{\delta_n+\sum_{k=1}^d|z_k|^{\alpha_k}}\d
 z.
 \eeqnn
 Since
$n\delta_n\to\theta\in[0, \infty)$ and $F_n^2=n$, from the above
formula, it follows that
 \beqlb\label{lem-3-7}
 \lim_{n\to\infty}\tilde{I}_{21}(n,\psi_n)=\frac{1}{(2\pi)^{d}}\int_0^1h(t)\d t\int_0^1 h(t')\d
 t'\int_0^{t\wedge t'}\e^{-\theta s}\d s
 \int_{\R^d}\frac{|\widehat{\phi}(z)|^2}{\sum_{k=1}^d|z_k|^{\alpha_k}}\d
 z.\qquad
 \eeqlb
To determine the limit of $I_{22}(n, \psi_n)$ in (\ref{lem-3-4}), we note the similarity
between $I_{22}(n, \psi_n)$ and $I_{12}(n, \psi_n)$. By using similar (but much simpler)
arguments, we can show
 \beqlb\label{lem-3-8}
 \lim_{n\to\infty}I_{22}(n,\psi_n)=0;
 \eeqlb
the details are omitted here. Combining (\ref{lem-3-4})
with (\ref{lem-3-6})-(\ref{lem-3-8}) yields (\ref{lem-3-1}).

To prove Part (2), let
\beqlb\label{lem-3-2}
\tilde{I}_2(n, \psi_n):=\int_{\R^d}\d x\int_0^n\e^{-\delta_n s}\psi_n(x,
 s) J_{n,\psi_n}(x, n-s, s)\d s.
 \eeqlb
Then (\ref{s3-1-13}) and (\ref{s3-1-18}) imply that
 \beqnn
 I_2(n,\psi_n)&\leq&\tilde{I}_2(n,\psi_n)=\int_0^n\e^{-\delta_n s}\d s
 \int_0^{n-s}f_n(v)\d v\int_{\R^d}\psi_n(x,
 s)T_v\psi_n(x, s+v)\d x.\qquad
 \eeqnn
By (\ref{s2-4}), (\ref{s3-2}), (\ref{s2-12}) and (\ref{s2-14}), there
exists a constant $K>0$ such that
\beqlb\label{lem-3-3}
 I_2(n,\psi_n)&\leq&\frac{K}{F_n^2}\int_0^n\e^{-\delta_n s}\d s\int_0^{n-s}
 \e^{-\delta_n v}\d v\int_{\R^d}|\widehat{\phi}(z)|^2\e^{-v\sum_{k=1}^d|z_k|^{\alpha_k}}\d
 z\nonumber
 \\&=&\frac{K}{F_n^2}\int_{\R^d}|\widehat{\phi}(z)|^2\d z\int_0^1\e^{-n\delta_n s}n
 \d s\int_0^{1-s}\e^{-n\delta_n v}n\e^{-nv\sum_{k=1}^d|z_k|^{\alpha_k}}\d
 v\nonumber
 \\&\leq&\frac{nK}{F_n^2}\int_{\R^d}\frac{|\widehat{\phi}(z)|^2}
 {\sum_{k=1}^d|z_k|^{\alpha_k}}\d z\int_0^1\e^{-n\delta_n s}\d s.
 \eeqlb
By Remark \ref{s2-rm-1},
$\int_{\R^d}\frac{|\widehat{\phi}(z)|^2}{\sum_{k=1}^d|z_k|^{\alpha_k}}\d
z<\infty$ for $\bar{\alpha}> 1$. Substituting
$F_n^2=n^{3-\bar{\alpha}}$ as $\bar{\alpha}\in (1, 2)$, or $F_n^2=n\ln n$
as $\bar{\alpha}=2$ into (\ref{lem-3-3}), we can readily see that
$  I_2(n,\psi_n)\to 0$ as $n\to\infty$. This finishes the proof.
\qed

The last lemma is concerned with $I_{3}(n,\psi_n)$ in (\ref{s3-1-19})
\begin{lem}\label{s3-lem-4}
If $F_n\to\infty$ as $n\to\infty$, then
 \beqlb\label{lem-4-1}
 \lim_{n\to\infty}I_{3}(n,\psi_n)=0.
 \eeqlb
\end{lem}
\medskip
{\bf Proof.}  From (\ref{s3-1-19}) and (\ref{s3-1-15}), we can obtain that
 \beqlb\label{lem-4-2}
 I_3(n,\psi_n)&\leq&\tilde{I}_3(n,\psi_n),
 \eeqlb
 where $\tilde{I}_3(n,\psi_n)$ is defined as
  \beqlb\label{lem-4-3}
 &&\delta_n\int_{\R^d}\d x\int_0^n\e^{-\delta_n s}\d
 s\int_0^{n-s}\e^{-\gamma u}\bfE_x\Big(\int_0^u\psi_n(\vec{\xi}_n(v), s+v)\d v
 \psi_n(\vec{\xi}_n(u), s+u)\Big)\d u\nonumber
 \\&&=\delta_n\int_0^n\e^{-\delta_n s}\d
 s\int_0^{n-s}\e^{-\gamma u}\d u\int_0^u\d v\int_{\R^d}\psi_n(x, s+v)T_{u-v}\psi_n(x,
 s+u)\d x.
 \eeqlb
By using (\ref{s2-12}), (\ref{s2-14}) and (\ref{s3-2}), we get from
(\ref{lem-4-3}) that
 \beqlb\label{lem-4-4}
  \tilde{I}_3(n,\psi_n)&\leq&\frac{\delta_n}{F_n^2(2\pi)^{d}}
  \int_0^n\e^{-\delta_n s}\d s\int_0^{n-s}\e^{-\gamma u}\d u
  \int_0^u \tilde{h}\Big(\frac{s+v}{n}\Big)\tilde{h}\Big(\frac{s+u}{n}\Big)\d
  v\nonumber
  \\&&\times\int_{\R^d}|\widehat{\phi}(z)|^2
  \e^{-(u-v)\sum_{k=1}^d|z_k|^{\alpha_k}}\d z.
 \eeqlb
The boundedness of $|\widehat{\phi}(z)|$ and $\tilde{h}$ implies that
 \beqlb\label{lem-4-5}
 \tilde{I}_3(n,\psi_n)\leq K\frac{\delta_n}{F_n^2}\int_0^n\e^{-\delta_n s}
 \d s\int_0^{n-s}\e^{-\gamma u}u\d u\leq\frac{K}{\gamma F_n^2}.
 \eeqlb
Therefore, (\ref{lem-4-1}) follows from (\ref{lem-4-2}) and (\ref{lem-4-5}).\qed

\section{Proofs of the main results}

In this section, we give the proofs of the main results stated in Section 2.

{\bf Proof of Theorem 2.1}. Without loss of generality, we prove the conclusion
for $t=1$, namely, we prove that $\langle\tilde{X}_n, \psi\rangle$
converges in distribution to $\langle\tilde{X}, \psi\rangle$ for
all $\psi\in\mathcal{S}(\R^{d+1})$, where $\tilde{X}_n$ and $\tilde{X}$ are
defined as in (\ref{s2-15}).
As explained in Bojdecki et al. \cite[p.9]{BGT061}, it suffices to show that
 \beqlb\label{s4-1}
 \lim_{n\to\infty}\bfE(\e^{-\langle\tilde{X}_n, \psi
 \rangle})=\exp\bigg(\frac{1}{2}\int_0^1\int_0^1 {\rm Cov}
 \Big(\big\langle X(s), \psi(\cdot, s)\big\rangle,\big\langle X(t),
 \psi(\cdot, t)\big\rangle\Big)\d s\d
 t\bigg)
 \eeqlb
for every non-negative $\psi\in\mathcal{S}(\R^{d+1})$. We only consider
the case of $\psi(x, t)=\phi(x)h(t)$.

 It follows from (\ref{s3-1-16}), (\ref{s3-2-1})--(\ref{s3-2-3})
 and Lemma \ref{s3-lem-1}-\ref{s3-lem-4}  that
 \beqnn
 \lim_{n\to\infty}\bfE(\e^{-\langle\tilde{X}_n, \psi
 \rangle})&=&\exp\bigg\{\frac{1}{(2\pi)^{d}}\int_0^1h(t)\d t\int_0^1 h(t')\d
 t'\int_0^{t\wedge t'}\e^{-\theta s}\d s
 \\&&\qquad\quad\times \int_{\R^d}\Big(\frac{1}{\sum_{k=1}^d|z_k|^{\alpha_k}} +\frac{\gamma/2}{(\sum_{k=1}^d|z_k|^{\alpha_k})^2}\Big)|\widehat{\phi}(z)|^2\d
 z\bigg\}.
 \eeqnn
Note that for the $\mathcal{S}'(\R^d)$-valued process $X$ with covariance (\ref{s2-thm-1-1}), we have
 \beqnn
 &&{\rm Cov}\Big(\big\langle X(s), \psi(\cdot, s)\big\rangle,\big\langle X(t),\psi(\cdot, t)\big\rangle\Big)
 \\&&\qquad=\frac{1}{(2\pi)^{d}}\bigg[\int_{\R^d}\bigg(\frac{2}{\sum_{k=1}^d|z_k|^{\alpha_k}}
 +\frac{\gamma}{(\sum_{k=1}^d|z_k|^{\alpha_k})^2}\bigg)\widehat{\phi}_1(z)\overline{\widehat{\phi}_2(z)}\d
 z\int_0^{r\wedge t}\e^{-\theta s}\d s\bigg] h(t)h(r).
 \eeqnn
Therefore, (\ref{s4-1}) holds.

For general $\psi\in \mathcal{S}(\R^{d+1})$, the proof is the same
with slightly more complicated notation. The details are omitted and
hence the proof of Theorem 2.1 is complete.\qed

{\bf Proof of Theorem \ref{s2-thm-2}}. The idea is same as that of Theorem \ref{s2-thm-1}.
The details are omitted.\qed

Below we prove Theorem \ref{s2-thm-3}.

{\bf Proof of Theorem \ref{s2-thm-3}}. We employ the space-time method formulated
in Bojdecki et al. \cite{BGT86}. Following Bojdecki et
al. \cite{BGT062}, it suffices to show the following two claims.

(i) $\langle\tilde{X}_n, \psi\rangle$ converges in distribution to
 $\langle \tilde{X},\psi\rangle$ for all $\psi\in\mathcal{S}(\R^{d+1})$ as
 $n\to\infty$.

(ii) $\{\langle X_n, \phi\rangle; n\geq 1\}$ is tight in $C([0, 1], \R)$
for all $\phi\in\mathcal{S}(\R^{d})$, where the
theorem of Mitoma \cite{M83} is used.

The proof of (i) is similar to that of Theorem \ref{s2-thm-1}. We sketch it
briefly as follows.
By  (\ref{s3-1-16}), (\ref{s3-2-1})-(\ref{s3-2-3}) and Lemma \ref{s3-lem-1}--\ref{s3-lem-4}
we can readily get
 \beqnn
 \lim_{n\to\infty}\bfE(\e^{-\langle\tilde{X}_n, \psi
 \rangle})&=&\exp\bigg\{\frac{\gamma}{2\pi^{d}}\prod_{k=1}^d
 \frac{\Gamma(1/\alpha_k)}{\alpha_k}\bigg(\int_{\R^d}\phi(x)\d
 x\bigg)^2\nonumber
 \\&&\qquad \times\int_0^1\e^{-\theta u}
 \tilde{h}(u)\d u\int_0^1\e^{-\theta v}\tilde{h}(v)\d v\int_0^{u\wedge v}
 \frac{\e^{\theta s}\d s}{(u+v-2s)^{\bar{\alpha}}}\bigg\}\nonumber
 \\&=&\exp\bigg\{\frac{\gamma}{2\pi^{d}}\prod_{k=1}^d\frac{\Gamma(1/\alpha_k)}
 {\alpha_k}\Big(\int_{\R^d}\phi(x)\d
 x\Big)^2\int_0^1\int_0^1 h(s)h(t)C(s, t)\d s\d t\bigg\}.\nonumber
 \eeqnn
The last term is the right hand side of (\ref{s4-1}) for the process $X$
in Theorem \ref{s2-thm-3} and $\psi(x,t)=\phi(x)h(t)$. This verifies (i).

Next we prove (ii). By Theorem 3.1 of Mimato \cite{M83} and the same argument
as that used in the proof of \cite[Proposition 3.3]{BGT071}, it suffices to prove that
for all $\phi\in\mathcal{S}(\R^d)$, $0\leq t_1<t_2\leq 1$ and
$\eta>0$, there exist constants $a\geq1$, $b>0$ and $K>0$, which is
independent of $t_1, t_2$, such that for all $n\geq1$.
 \beqlb\label{s4-2}
 \int_0^{1/\eta}\Big(1-{\rm Re}\Big(\bfE\Big[\exp\{-i\omega\langle\tilde{X}_n,
 \phi h\rangle\}\Big]\Big)\Big)\d\omega\leq\frac{K}{\eta^{a}}(t_2-t_1)^{1+b},
 \eeqlb
where $h\in\mathcal{S}(R)$ is an approximation of ${\bf
1}_{\{t_2\}}(t)-{\bf 1}_{\{t_1\}}(t)$ supported on $[t_1, t_2]$ such
that $\tilde{h}(t)$ satisfies
 \beqlb\label{s4-3}
 \tilde{h}\in\mathcal{S}(R),\qquad 0\leq \tilde{h}\leq {\bf 1 }_{[t_1, t_2]}.
 \eeqlb
We now repeat the argument in Section 3 with $\phi$ replaced by
$i\omega\phi, \omega>0$ and $h$ satisfying (\ref{s4-3}) and derive
$$
\bfE\Big[\exp\Big\{-i\omega\langle\tilde{X}_n,
\phi h\rangle\Big\}\Big]=\exp\Big\{I_1(n, i\omega\psi_n)+
I_2(n, i\omega\psi_n)+I_3(n, i\omega\psi_n)\Big\}
 $$
and
$$
|V_{n, i\omega\psi_n}|\leq J_{n,\omega\psi_n}.
$$
Consequently, from the expressions of $I_1, I_2, I_3$ and $I_{11}$,
(\ref{lem-3-2}) and (\ref{lem-4-2}),  we can verify that
\begin{equation}\label{s4-4}
\begin{split}
 &|I_1(n, i\omega\psi_n)|\leq I_{11}(n,
 \omega\psi_n); \quad |I_2(n, i\omega\psi_n)|\leq \tilde{I}_2(n,
 \omega\psi_n);  \\
 & |I_3(n, i\omega\psi_n)|\leq \tilde{I}_3(n,  \omega\psi_n).
 \end{split}
 \end{equation}
In the following we estimate $I_{11}(n, \omega\psi_n)$, $\tilde{I}_2(n,
\omega\psi_n)$ and $ \tilde{I}_3(n,  \omega\psi_n)$ separately.

For $I_{11}(n, \omega\psi_n)$, (\ref{s3-2}), (\ref{s3-1-13}) and the
boundedness of $\{f_n\}_{n\geq 1}$ imply that for some constant $K>0$,
\beqnn
 I_{11}(n,\omega\psi_n)&\leq& \frac{K\omega^2}{F_n^2}\int_0^n\d
 s\int_0^{n-s}\tilde{h}\Big(\frac{s+u}{n}\Big)\d
 u\int_0^{n-s}\tilde{h}\Big(\frac{s+v}{n}\Big)\d v
 \int_{\R^d} T_u\phi(x)T_v\phi(x)\d x\nonumber
 \\&=&\frac{2K\omega^2n^3}{F_n^2}\int_0^1\d
 s\int_s^1\tilde{h}(u)\d
 u\int_u^1\tilde{h}(v)\d v\int_{\R^d} |\widehat{\phi}(z)|^2\e^{-n(u+v-2s)\sum_{k=1}^d|z_k|^{\alpha_k}}\d z
 \\&=&\frac{2K\omega^2n^3}{F_n^2}\int_{\R^d} |\widehat{\phi}(z)|^2\d z\int_0^1\tilde{h}(u)\d
 u\int_0^u\tilde{h}(v)\d v\int_0^v\e^{-n(u+v-2s)\sum_{k=1}^d|z_k|^{\alpha_k}}\d
 s,
 \eeqnn
which, combined with (\ref{s4-3}), implies that $I_{11}(n,\psi_n)$
is bounded from above by
\beqlb\label{s4-5}
 &&\frac{2K\omega^2n^3}{F_n^2}\int_{\R^d} |\widehat{\phi}(z)|^2\d
 z\int_{t_1}^{t_2}\d
 u\int_{t_1}^u\e^{-n(u+v)\sum_{k=1}^d|z_k|^{\alpha_k}}\d v\int_0^v\e^{2ns\sum_{k=1}^d|z_k|^{\alpha_k}}\d
 s\nonumber
 \\&&\qquad\leq\frac{K\omega^2n^2}{F_n^2}\int_{\R^d} \frac{|\widehat{\phi}(z)|^2}{\sum_{k=1}^d|z_k|^{\alpha_k}}\d
 z\int_{t_1}^{t_2}\d
 u\int_{t_1}^u\e^{-n(u-v)\sum_{k=1}^d|z_k|^{\alpha_k}}\d v \nonumber
 \\&&\qquad=\frac{K\omega^2n}{F_n^2}\int_{\R^d} \frac{|\widehat{\phi}(z)|^2}{\sum_{k=1}^d|z_k|^{\alpha_k}}\d
 z\int_{t_1}^{t_2}\frac{1-\e^{-n(u-t_1)\sum_{k=1}^d|z_k|^{\alpha_k}}}{\sum_{k=1}^d|z_k|^{\alpha_k}}\d
 u.
\eeqlb
Substituting $F_n^2=n^{3-\bar{\alpha}}$ and $y=n^Hz$  into (\ref{s4-5}), we get
\beqnn
 I_{11}(n,\omega\psi_n)&\leq&K\omega^2\int_{\R^d} \frac{|\widehat{\phi}(n^{-H}y)|^2}{\sum_{k=1}^d|y_k|^{\alpha_k}}\d
 y\int_{t_1}^{t_2}\frac{1-\e^{-(u-t_1)\sum_{k=1}^d|y_k|^{\alpha_k}}}{\sum_{k=1}^d|y_k|^{\alpha_k}}\d
 u.\nonumber
 \eeqnn
By using the inequality $1-\e^{-x}\leq x^r$ for all $x\geq0$ and $r\in(0, 1]$ and the boundedness of
$|\widehat{\phi}|$,  we obtain that,
\beqnn
 I_{11}(n,\omega\psi_n)&\leq&K\omega^2 \|\widehat{\phi}\|^2\bigg[\int_{[0, 1]^d}
 \frac{1}{(\sum_{k=1}^d|z_k|^{\alpha_k})^{2-r_1}}\d
 z\int_{t_1}^{t_2}(u-t_1)^{r_1}\d u\nonumber
 \\&&+K\int_{\R^d\setminus[0, 1]^d} \frac{1}{(\sum_{k=1}^d|z_k|^{\alpha_k})^{2-r}}\d
 z\int_{t_1}^{t_2}(u-t_1)^{r}\d u\bigg],
 \eeqnn
where $\|\widehat{\phi}\|:=\sup_{z\in\R^d}|\widehat{\phi}(z)|$. Since $\bar{\alpha}\in(1,
 2)$, Lemma \ref{s2-lm-1} and  Remark \ref{s2-rm-1} imply that for any
 $r_1\in(2-\bar{\alpha}, 1)$ and $r\in(0, 2-\bar{\alpha})$,
  $$\int_{[0, 1]^d} \frac{1}{(\sum_{k=1}^d|z_k|^{\alpha_k})^{2-r_1}}\d
 z+\int_{\R^d\setminus[0, 1]^d} \frac{1}{(\sum_{k=1}^d|z_k|^{\alpha_k})^{2-r}}\d
 z<\infty.$$
Therefore, there exist a constant $K>0$ such that
\beqlb\label{s4-6}
I_{11}(n,\omega\psi_n)&\leq&
K\omega^2\|\widehat{\phi}\|^2|t_2-t_1|^{1+r}.
\eeqlb

Next we estimate $\tilde{I}_2(n,\omega\psi_n)$. Due to the boundedness of $\bar{f}_n$,
(\ref{lem-3-3}) implies that for some
constant $K>0$
 \beqnn
 \tilde{I}_2(n,\omega\psi_n)&\leq& \frac{K\omega^2}{F_n^2}\int_0^n\tilde{h}\Big(\frac{s}{n}\Big)\d
 s\int_0^{n-s}\tilde{h}\Big(\frac{s+u}{n}\Big)\d
 u\int_{\R^d} \phi(x)T_u\phi(x)\d x.
 \eeqnn
Then by condition (\ref{s4-3}), we have
 \beqlb\label{s4-7}
 \tilde{I}_2(n,\omega\psi_n)&\leq&\frac{K\omega^2n^2}{F_n^2}\int_{t_1}^{t_2}\d
 s\int_{s}^{t_2}\d
 u\int_{\R^d} |\widehat{\phi}(z)|^2\e^{-n(u-s)\sum_{k=1}^d|z_k|^{\alpha_k}}\d
 z\nonumber
  \\&\leq&\frac{K\omega^2n}{F_n^2}\int_{\R^d} \frac{|\widehat{\phi}(z)|^2}{\sum_{k=1}^d|z_k|^{\alpha_k}}\d
 z\int_{t_1}^{t_2}(1-\e^{-n(t_2-s)\sum_{k=1}^d|z_k|^{\alpha_k}})\d
 s.
 \eeqlb
Using the inequality $1-\e^{-x}\leq x^r$ for all $x\geq0$ and $r\in(0, 1]$ again,
 and substituting $F_n^2=n^{3-\bar{\alpha}}$ into (\ref{s4-7}), we get that,
 \beqnn
 \tilde{I}_2(n,\omega\psi_n)&\leq&\frac{K\omega^2n^{1+r}}{n^{3-\bar{\alpha}}} \int_{\R^d}
 \frac{|\widehat{\phi}(z)|^2}{(\sum_{k=1}^d|z_k|^{\alpha_k})^{1-r}}\d
 z\int_{t_1}^{t_2}(u-t_1)^r\d u.
 \eeqnn
Let $r=2-\bar{\alpha}$. Then $1-r=\bar{\alpha}-1\in(0, 1)$ and $\int_{\R^d}
\frac{|\widehat{\phi}(z)|^2}{(\sum_{k=1}^d|z_k|^{\alpha_k})^{1-r}}\d
 z<\infty$, thanks to Remark \ref{s2-rm-1}. Therefore there exists a
 constant $K>0$ such that
 \beqlb\label{s4-8}
 \tilde{I}_2(n,\omega\psi_n)&\leq&K\omega^2\int_{\R^d}
 \frac{|\widehat{\phi}(z)|^2}{(\sum_{k=1}^d|z_k|^{\alpha_k})^{1-r}}\d
 z |t_2-t_1|^{1+r}.
 \eeqlb

In order to estimate $\tilde{I}_3(n,\omega\psi_n)$, we combine (\ref{lem-4-4})
and (\ref{s4-3}) to see that for all $r\in(0, 1)$,
  \beqnn
  \tilde{I}_3(n,\omega\psi_n)&\leq&\frac{\omega^2\delta_n}{F_n^2(2\pi)^{d}}\int_0^n\d s\int_0^{n-s}\e^{-\gamma u}\d u
  \int_0^u \tilde{h}\Big(\frac{s+v}{n}\Big)\tilde{h}\Big(\frac{s+u}{n}\Big)\d v\int_{\R^d}|\widehat{\phi}(z)|^2\d z\nonumber
  \\&=&\frac{\omega^2n^3\delta_n}{F_n^2(2\pi)^{d}}\int_{\R^d}|\widehat{\phi}(z)|^2\d z\int_0^1
  \d s\int_s^{1}\e^{-n\gamma(u-s)}\tilde{h}(u)\d u
  \int_s^u \tilde{h}(v)\d v\nonumber
  \\&\leq&\frac{\omega^2n^3\delta_n}{F_n^2}K\int_{t_1}^{t_2}\e^{-n\gamma u}\d u\int_{t_1}^u\d v
  \int_0^v\e^{ns\gamma}\d s\nonumber
  \\&\leq&\frac{\omega^2n\delta_n}{\gamma^2F_n^2}K
  \int_{t_1}^{t_2}(1-\e^{-n\gamma(u-t_1)})\d
  u\\
  &\leq&\frac{\omega^2\gamma^rn^{1+r}\delta_n}{\gamma^2F_n^2}K|t_2-t_1|^{1+r}.
 \eeqnn
From $n\delta_n\to\theta\in[0, \infty)$, $F_n^2=n^{3-\bar{\alpha}}$ and $1 <\bar{\alpha}<2$
we derive that
  $$0<\frac{\gamma^rn^{1+r}\delta_n}{\gamma^2F_n^2}\leq
  \frac{\gamma^rn\delta_n}{\gamma^2n^{3-r-\bar{\alpha}}}\to 0.$$
Therefore,  for any $r\in(0,  1)$ there exists a constant $K$ such that
 \beqlb\label{s4-9}
 \tilde{I}_3(n,\omega\psi_n)\leq K\omega^2|t_2-t_1|^{1+r}.
 \eeqlb

Combining (\ref{s4-6}) with (\ref{s4-8}) and (\ref{s4-9}), we have that for some
$r\in(0, 2-\bar{\alpha})$, there is a constant $K$
independent of $t_1, t_2$ and $r>0$ such that
 \beqlb\label{s4-10}
 |\tilde{I}_3(n, \omega\psi_n)|+|\tilde{I}_2(n, \omega\psi_n)|+|I_{11}
 (n, \omega\psi_n)|\leq K(\phi, r)\omega^2|t_2-t_1|^{1+r}.
 \eeqlb
Note that
$$
\Big|1-{\rm Re}\Big(\bfE\Big[\exp\big\{-i\omega\langle\tilde{X}_n,
 \phi h\rangle\big\}\Big]\Big)\Big|\leq |I_1(n, i\omega\psi_n)|+
 |I_2(n, i\omega\psi_n)|+|I_3(n, i\omega\psi_n)|,$$
we derive from (\ref{s4-4}) and (\ref{s4-10}) that
$$
\int_0^{1/\eta}\Big(1-{\rm Re}\Big(\bfE\Big[\exp\big\{-i\omega\langle\tilde{X}_n,
 \phi h\rangle\big\}\Big]\Big)\Big)\d\omega\leq\frac{K(\phi,
 r)}{3\eta^3}|t_2-t_1|^{1+r}.$$
This completes the proof of (\ref{s4-2}) and hence the
proof of Theorem \ref{s2-thm-1}.\qed

\medskip

At last, we prove Proposition \ref{s2-prop}.

{\bf Proof of Proposition \ref{s2-prop}.}  We only prove the statements on $Y_1$.
The remainder is similar and omitted.

By the definition of $Y_1 = \{Y_1(u), u \in [0, \infty)^d\}$, it is a centered
Gaussian random field with covariance function given by
 \beqlb\label{Eq:Y1}
 \bfE[Y_1(u)Y_1(v)]&=&\frac{2}{(2\pi)^d}\int_0^1\e^{-\theta s}\d s
 \int_{\R^d} \frac{\d z} {\sum_{k=1}^d|z_k|^{\alpha_k}}
 \int_{D(u)}\int_{D(v)} \e^{-i\sum_{k=1}^d(x_k-y_k)z_k} \d x\d y \nonumber
 \\&=&K\,
 \int_{\R^d} \prod_{k=1}^d\frac{(1-\e^{iu_kz_k})(1-\e^{-iv_kz_k})}{z_k^2}
 \frac{\d z}{\sum_{k=1}^d|z_k|^{\alpha_k}}.
 \eeqlb
Here and below, $K=2(1-\e^{-\theta})/((2\pi)^d\theta)$. Note that the last integral is
finite because $\bar{\alpha}>2$.

{\bf (1)} If $d=1$, then $\vec{\alpha}=(\alpha_1)=:\alpha<1/2$.
It follows from (\ref{Eq:Y1}) that for any $u,v\in\R_+$,
 \beqnn
 \bfE[Y_1(u)Y_1(v)]&=&K \int_{\R} \frac{(1-\e^{iuz})(1-\e^{-ivz})}{|z|^{2+\alpha}}
 \d z
 \\&=& K  \int_{\R} \frac{1-\cos (uz)-\cos (vz)+\cos((u-v)z)}{|z|^{2+\alpha}}\d z
 \\&=&K  \Big[u^{1+\alpha}+v^{1+\alpha}-|u-v|^{1+\alpha}\Big]\int_{\R}\frac{1-\cos x}{|x|^{2+\alpha}}\d x.
  \eeqnn
It is well-known that the fractional Brownian motion $B_h$ with Hurst exponent
$h\in(0, 1)$ is a one-parameter centered Gaussian process  with covariance
 $$\bfE(B_h(u)B_h(v))=\frac{1}{2}\Big(u^{2h}+v^{2h}-|u-v|^{2h}\Big).$$
Therefore, up to a multiplicative constant, $Y_1$ is the fractional Brownian
motion with Hurst exponent $(1+\alpha)/2\in(1/2, 3/4)$.

{\bf (2)} Now assume that $d > 1$. Recall that $H$ is the $d\times d$ diagonal
matrix $(1/\alpha_k)_{1\leq k\leq d}$. For any $r > 0$, it follows from (\ref{Eq:Y1})
and the substitution $x=r^Hz$ that
\beqnn
\bfE[Y_1(r^Hu)Y_1(r^Hv)]&=&K \int_{\R^d} \prod_{k=1}^d
\frac{(1-\e^{ir^{1/\alpha_k}u_kz_k})(1-\e^{-ir^{1/\alpha_k}v_kz_k})}{z_k^2}
\frac{\d z}{\sum_{k=1}^d|z_k|^{\alpha_k}}\\
&=& r^{1+\bar{\alpha}}K\,\int_{\R^d} \prod_{k=1}^d\frac{(1-\e^{iu_kx_k})(1-\e^{-iv_kx_k})}{x_k^2}
\frac{\d x}{\sum_{k=1}^d|x_k|^{\alpha_k}}
\\&=&r^{1+\bar{\alpha}}\bfE[Y_1(u)Y_1(v)].
  \eeqnn
Therefore, for any $r > 0$,
$$
\{Y_1(r^H u), u\in[0,\infty)^d\}\stackrel{f.d}{=}
\{r^{(\bar{\alpha}+1)/2}Y_1(u), u\in[0,\infty)^d\},
$$
That is, $Y_1$ is an operator scaling random field with exponent $H$.

Let $v_1=(0,0,\cdots,0)$, $v_2=(1,0,\cdots,0)$ and $u_1=(0,1,\cdots,1)$, $u_2=(1,1,\cdots,1)$.
Then $u_1 - v_1 = u_2 - v_2$ and (\ref{Eq:Y1}) gives
$$ \bfE[(Y_1(v_1))^2]= \bfE[(Y_1(u_1))^2]= \bfE[(Y_1(v_2))^2]= \bfE[Y_1(u_2)Y_1(v_2)]=0,$$
Hence
\beqnn
 &&\bfE[(Y_1(u_2)-Y_1(v_2))^2]=\bfE[(Y_1(u_2))^2]
 \\&&\qquad=K \, \int_{\R^d} \prod_{k=1}^d\frac{(1-\e^{ix_k})(1-\e^{-ix_k})}{x_k^2}
 \frac{\d x}{\sum_{k=1}^d|x_k|^{\alpha_k}}
  \\&&\qquad\not=0=\bfE[(Y_1(u_1))^2]=\bfE[(Y_1(u_1)-Y_1(v_1))^2],
  \eeqnn
which implies that $Y_1$ does not have stationary increments. \qed

\bigskip

\ack
The authors are grateful to the anonymous referees and Associate Editor
for careful reading of the paper and for their
suggestions which have helped to improve the paper significantly.

\end{document}